\documentclass[a4paper, 11pt]{article}

\usepackage{xspace}
\usepackage{tikz}
\usepackage{subfigure}
\usepackage{rotating}
\usepackage{amsmath}
\usepackage{amsfonts}
\usepackage{url}
\usepackage{enumitem}

\def\bfa{\mbox{\boldmath$a$}}

\def\bfc{\mbox{\boldmath$c$}}

\def\bfe{\mbox{\boldmath$e$}}

\def\bfl{\mbox{\boldmath$l$}}

\def\bfu{\mbox{\boldmath$u$}}

\def\bfx{\mbox{\boldmath$x$}}

\def\bfzero{\mbox{\boldmath$0$}}

\def\bftau{\mbox{\boldmath$\tau$}}

\def\SetB{{\cal B}}

\def\SetF{{\cal F}}

\def\SetJ{{\cal J}}

\def\SetN{{\cal N}}
\def\SetP{{\cal P}}

\def\subB{_{\scriptscriptstyle B}}
\def\subF{_{\scriptscriptstyle F}}
\def\subN{_{\scriptscriptstyle N}}

\newcommand{\hsol}{\texttt{hsol}\xspace}
\newcommand{\sip}{\texttt{sip}\xspace}
\newcommand{\sipeight}{\texttt{sip8}\xspace}
\newcommand{\pami}{\texttt{pami}\xspace}
\newcommand{\pamione}{\texttt{pami1}\xspace}
\newcommand{\pamieight}{\texttt{pami8}\xspace}
\newcommand{\clp}{\texttt{Clp}\xspace}
\newcommand{\cplex}{\texttt{Cplex}\xspace}
\newcommand{\xpress}{\texttt{Xpress}\xspace}
\newcommand{\xpresseight}{\texttt{Xpress8}\xspace}
\newcommand{\LP}[1]{\textsc{#1}\xspace}

\def\CHUZC{{\sc chuzc}\xspace}
\def\CHUZCONE{{\sc chuzc1}\xspace}
\def\CHUZCTWO{{\sc chuzc2}\xspace}
\def\CHUZR{{\sc chuzr}\xspace}
\def\FTRAN{{\sc ftran}\xspace}
\def\BTRAN{{\sc btran}\xspace}
\def\PRICE{{\sc spmv}\xspace}
\def\INVERT{{\sc invert}\xspace}
\def\UPDATE{{\sc update}\xspace}

\def\UPDATEDUAL{{\sc update-dual}\xspace}
\def\UPDATEPRIMAL{{\sc update-primal}\xspace}
\def\UPDATEWEIGHT{{\sc update-weight}\xspace}
\def\UPDATEFACTOR{{\sc update-factor}\xspace}
\def\FTRANDSE{{\sc ftran-dse}\xspace}
\def\FTRANBFRT{{\sc ftran-bfrt}\xspace}
\def\OTHER{{\sc other}\xspace}

\def\LP#1{\textsc{#1}\xspace}

\newcommand{\inred}[1]{{\color{red} #1}}

\def\nullbox#1#2{\setbox0=\null\ht0=#1 \dp0=#2\box0}

\tikzstyle{block} = [rectangle, draw, fill=gray!20, font=\small\ttfamily,
    text width=2cm, text centered, minimum height=1cm]

\newcommand{\TableOneForceLineHeight}{\nullbox{12pt}{2pt}}

\newcommand{\JHchange}[2]{{\textsl{}}{\inred{}}}

\newcommand{\JHcomment}[2]{{\textsl{}}{\inred{}}}

\title{Parallelizing the dual revised simplex method}
\author{
Q.~Huangfu\thanks{FICO House, International Square, Starley Way, Birmingham, B37 7GN, UK} and 
J.~A.~J.~Hall\thanks{School of Mathematics and Maxwell Institute for Mathematical Sciences, University of Edinburgh, James Clerk Maxwell Building, Peter Guthrie Tait Road, Edinburgh, EH9 3FD, UK. Tel.: [+44](131) 650 5075, Fax: [+44](131) 650 6553; Email J.A.J.Hall@ed.ac.uk}}

\newcommand{\tbl}[2]{\caption{#1}#2}

\begin{document}
\maketitle

\begin{abstract}
This paper introduces the design and implementation of two parallel dual simplex solvers for general large scale sparse linear programming problems.  One approach, called PAMI, extends a relatively unknown pivoting strategy called suboptimization and exploits parallelism across multiple iterations.  The other, called SIP, exploits purely single iteration parallelism by overlapping computational components when possible.  Computational results show that the performance of PAMI is superior to that of the leading open-source simplex solver, and that SIP complements PAMI in achieving speedup when PAMI results in slowdown. One of the authors has implemented the techniques underlying PAMI within the FICO Xpress simplex solver and this paper presents computational results demonstrating their value. This performance increase is sufficiently valuable for the achievement to be used as the basis of promotional material by FICO. In developing the first parallel revised simplex solver of general utility and commercial importance, this work represents a significant achievement in computational optimization.

\end{abstract}
\hspace{0.5cm}

\noindent {\it Keywords:}  
Revised simplex method,
simplex parallelization

\section{Introduction}

%
%
Linear programming (LP) has been used widely and successfully in many
practical areas since the introduction of the simplex method in the
1950s.  Although an alternative solution technique, the interior point
method (IPM), has become competitive and popular since the 1980s, the
dual revised simplex method is frequently preferred, particularly when
families of related problems are to be solved.

%
The standard simplex method implements the simplex algorithm via a
rectangular tableau but is very inefficient when applied to sparse LP
problems. For such problems the revised simplex method is preferred
since it permits the (hyper-)sparsity of the problem to be
exploited. This is achieved using techniques for factoring sparse
matrices and solving hyper-sparse linear systems. Also important for
the dual revised simplex method are advanced algorithmic variants
introduced in the 1990s, particularly dual steepest-edge (DSE) pricing
and the bound flipping ratio test (BFRT). These led to dramatic
performance improvements and are key reasons for the dual simplex algorithm
being preferred.

%
A review of past work on parallelising the simplex method is given by
Hall~\cite{Hall2010}. The standard simplex method has been
parallelised many times and generally achieves good speedup, with
factors ranging from tens to up to a thousand. However, without using
expensive parallel computing resources, its performance on sparse LP problems is inferior to
a good sequential implementation of the revised
simplex method. The standard simplex method is also unstable numerically. Parallelisation of the revised simplex method has
been considered relatively little and there has been less success in
terms of speedup. Indeed, since scalable speedup for general large
sparse LP problems appears unachievable, the revised simplex method
has been considered unsuitable for parallelisation. However, since
it corresponds to the computationally efficient serial technique, any
improvement in performance due to exploiting parallelism in the
revised simplex method is a worthwhile goal.

%
Two main factors motivated the work in this paper to develop a
parallelisation of the dual revised simplex method for standard desktop
architectures. Firstly, although dual simplex implementations are now
generally preferred, almost all the work by others on parallel simplex has been
restricted to the primal algorithm, the only published work on dual
simplex parallelisation known to the authors being due to Bixby and
Martin~\cite{Bixby2000}. Although it appeared in the early 2000s,
their implementation included neither the BFRT nor hyper-sparse linear
system solution techniques so there is immediate scope to extend their
work. Secondly, in the past, parallel implementations generally used dedicated
high performance computers to achieve the best performance. Now, when
every desktop computer is a multi-core machine, any speedup is
desirable in terms of solution time reduction for daily use. Thus we
have used a relatively standard architecture to perform computational
experiments.

%
A worthwhile simplex parallelisation should be based on a good
sequential simplex solver.  Although there are many public domain
simplex implementations, they are either too complicated to be used as
a foundation for a parallel solver or too inefficient for any
parallelisation to be worthwhile. Thus the authors have implemented a
sequential dual simplex solver (\hsol) from scratch. It incorporates
sparse LU factorization, hyper-sparse linear system solution
techniques, efficient approaches to updating LU factors and
sophisticated dual revised simplex pivoting rules. Based on components
of this sequential solver, two dual simplex parallel solvers (\pami
and \sip) have been designed and developed.
Section~\ref{sect:Background} introduces the necessary background,
Sections~\ref{sect:PAMI} and~\ref{sect:SIP} detail the design of \pami
and \sip respectively and Section~\ref{sect:Results} presents
numerical results and performance analysis.  Conclusions are given in
Section~\ref{sect:Conclusions}.

\section{Background}\label{sect:Background}

The simplex method has been under development for more than 60 years,
during which time many important algorithmic variants have enhanced
the performance of simplex implementations. As a result, for novel
computational developments to be of value they must be tested within
an efficient implementation or good reasons given why they are
applicable in such an environment. Any development which is only
effective in the context of an inefficient implementation is not
worthy of attention.

This section introduces all the necessary background knowledge for
developing the parallel dual simplex solvers.
Section~\ref{background:lp} introduces the computational form of LP
problems and the concept of primal and dual feasibility.
Section~\ref{background:dual} describes the regular dual simplex
method algorithm and then details its key enhancements and major
computational components.  Section~\ref{background:subopt} introduces
suboptimization, a relative unknown dual simplex variant which is the
starting point for the \pami parallelisation in
Section~\ref{sect:PAMI}. Section~\ref{background:update} briefly
reviews several existing simplex update approaches which are key to
the efficiency of the parallel schemes.

\subsection{Linear programming problems}
\label{background:lp}

%
%
A linear programming (LP) problem in general computational form is
\begin{equation}
\text{minimize } f = \bfc^T \bfx \quad
\text{ subject to } A \bfx = \bfzero
\text{ and } \bfl \leq \bfx \leq \bfu,
\label{equ:complp}
\end{equation}
where $A \in \mathbb{R}^{m\times n}$ is the coefficient matrix and
$\bfx$, $\bfc$, $\bfl$ and $\bfu\in\mathbb{R}^m$ are, respectively,
the variable vector, cost vector and (lower and upper) bound vectors.
Bounds on the constraints are incorporated into $\bfl$ and $\bfu$ via
an identity submatrix of $A$. Thus it may be assumed that $m < n$ and
that $A$ is of full rank.

%
%
As $A$ is of full rank, it is always possible to identify a
non-singular basis partition $B\in \mathbb{R}^{m\times m}$ 
consisting of $m$ linearly independent columns of $A$, with the
remaining columns of $A$ forming the matrix $N$.  The variables are
partitioned accordingly into basic variables $\bfx\subB$ and nonbasic
variables $\bfx\subN$, so $A \bfx = B\bfx\subB + N\bfx\subN =
\bfzero$, and the cost vector is partitioned into basic costs
$\bfc\subB$ and nonbasic costs $\bfc\subN$, so $f =
\bfc\subB^T\bfx\subB + \bfc\subN^T\bfx\subN$. The indices of the basic
and nonbasic variables form sets $\SetB$ and $\SetN$ respectively.

%
%
In the simplex algorithm, the values of the (primal) variables are
defined by setting each nonbasic variable to one of its finite bounds
and computing the values of the basic variables as $\bfx\subB = - B^{-1} N
\bfx\subN$. The values of the dual variables (reduced costs) are
defined as $\widehat{\bfc}\subN^T = \bfc\subN^T - \bfc\subB^T B^{-1}N$.
When $\bfl\subB \leq \bfx\subB \leq \bfu\subB$ holds, the basis is
said to be primal feasible.  Otherwise, the primal infeasibility for
each basic variable $i \in \SetB$ is defined as
\begin{equation}
\Delta x_i = \left\{ 
\begin{array} {ll}
l_i - x_i & \text{ if } x_i < l_i \\
x_i - u_i & \text{ if } x_i > u_i \\
0         & \text{ otherwise}
\end{array}
\right. 
\label{equ:primalinfeas}
\end{equation}
If the following condition holds for all $j\in\SetN$ such that $l_j \neq u_j$
\begin{equation}
\widehat{c}_j \geq 0~(x_j = l_j),\quad \widehat{c}_j \leq 0~(x_j = u_j)
\label{equ:dualfeas}
\end{equation}
then the basis is said to be dual feasible.  It can be proved that if
a basis is both primal and dual feasible then it yields an optimal
solution to the LP problem. 

\subsection{Dual revised simplex method}
\label{background:dual}

The dual simplex algorithm solves an LP problem iteratively by seeking
primal feasibility while maintaining dual feasibility.  Starting from
a dual feasible basis, each iteration of the dual simplex algorithm can
be summarised as three major operations.

\begin{enumerate}
  \item {\it Optimality test}. In a component known as \CHUZR, choose
    the index $p \in \SetB$ of a good primal infeasible variable
    to leave the basis. If no such variable can be chosen, the LP
    problem is solved to optimality.
  \item {\it Ratio test}.  In a component known as \CHUZC, choose the
    index $q \in \SetN$ of a good nonbasic variable to enter the basis
    so that, within the new partition, $\widehat{c}_q$ is zeroed whilst
    $\widehat{c}_p$ and other nonbasic variables remain dual feasible.
    This is achieved via a ratio test with $\widehat{\bfc}^T\subN$ and
    $\widehat{\bfa}_p^T$, where $\widehat{\bfa}_p^T$ is row $p$ of the reduced
    coefficient matrix $\widehat{A}=B^{-1}A$.
  \item {\it Updating}.  The basis is updated by interchanging indices
    $p$ and $q$ between sets $\SetB$ and $\SetN$, with corresponding
    updates of the values of the primal variables $\bfx\subB$ using
    $\widehat{\bfa}_q$ (being column $q$ of $\widehat{A}$) and dual
    variables~$\widehat{\bfc}^T\subN$ using~$\widehat{\bfa}_p^T$, as well
    as other components as discussed below.
\end{enumerate}

What defines the revised simplex method is a representation of the
basis inverse $B^{-1}$ to permit rows and columns of the reduced
coefficient matrix $\widehat{A}=B^{-1}A$ to be computed by solving
linear systems.  The operation to compute the representation of
$B^{-1}$ directly is referred to as \INVERT and is generally achieved
via sparsity-exploiting LU factorization.  At the end of each simplex
iteration the representation of $B^{-1}$ is updated until it is
computationally advantageous or numerically necessary to compute a
fresh representation directly.  The computational component which
performs the update of $B^{-1}$ is referred to as \UPDATEFACTOR.
Efficient approaches for updating $B^{-1}$ are summarised in
Section~\ref{background:update}.

For many sparse LP problems the matrix $B^{-1}$ is dense, so solutions
of linear systems involving $B$ or $B^T$ can be expected to be dense
even when, as is typically the case in the revised simplex method,
the RHS is sparse. However, for some classes of LP problem the
solutions of such systems are typically sparse. This phenomenon, and
techniques for exploiting in the simplex method, it was identified by
Hall and McKinnon~\cite{Hall2005} and is referred to as
hyper-sparsity. 

The remainder of this section introduces advanced
algorithmic components of the dual simplex method.

\subsubsection{Optimality test}
In the optimality test, a modern dual simplex implementation adopts two
important enhancements.  The first is the dual steepest-edge (DSE)
algorithm~\cite{Forrest1992} which chooses the basic variable with
greatest weighted infeasibility as the leaving variable. This variable
has index
\[
p = \arg \max_i \frac{\Delta x_i} {||\widehat{\bfe}_i^T||_2}.
\]
For each basic variable $i \in \SetB$, the associated DSE weight $w_i$
is defined as the 2-norm of row $i$ of $B^{-1}$ so $w_i =
||\widehat{\bfe}_i^T||_2 =||{\bfe}_i^TB^{-1}||_2$.  The weighted
infeasibility $\alpha_i = \Delta x_i / w_i$ is referred to as the
attractiveness of a basic variable.  The DSE weight is updated at the
end of the simplex iteration.

The second enhancement of the optimality test is the hyper-sparse
candidate selection technique originally proposed for column selection
in the primal simplex method~\cite{Hall2005}. This maintains a short
list of the most attractive variables and is more efficient for
large and sparse LP problems since it avoids repeatedly searching the less
attractive choices.  This technique has been adapted for the dual
simplex row selection component of \hsol.

\subsubsection{Ratio test}
In the ratio test, the updated pivotal row $\widehat{\bfa}_p^T$ is
obtained by computing $\widehat{\bfe}_p^T = \bfe_p^TB^{-1}$ and then
forming the matrix vector product $\widehat{\bfa}_p^T =
\widehat{\bfe}_p^T A$.  These two computational components are
referred to as \BTRAN and \PRICE respectively.

The dual ratio test (\CHUZC) is enhanced by the Harris two-pass ratio
test~\cite{Harris1973} and bound-flipping ratio test
(BFRT)~\cite{Fourer1994}.  Details of how to apply these two
techniques are set out by Koberstein~\cite{Koberstein2008}.

For the purpose of this report, advanced \CHUZC can be viewed as
having two stages, an initial stage \CHUZCONE which simply accumulates
all candidate nonbasic variables and then a recursive selection stage
\CHUZCTWO to choose the entering variable $q$ from within this set of
candidates using BFRT and the Harris two-pass ratio test. \CHUZC also
determines the primal step $\theta_p$ and dual step $\theta_q$, being
the changes to the primal basic variable $p$ and dual variable $q$
respectively.  Following a successful BFRT, \CHUZC also yields an
index set $\SetF$ of any primal variables which have flipped from one
bound to the other.

\subsubsection{Updating}
In the updating operation, besides \UPDATEFACTOR, several vectors are
updated. Update of the basic primal variables $\bfx\subB$
(\UPDATEPRIMAL) is achieved using $\theta_p$ and $\widehat{\bfa}_q$,
where $\widehat{\bfa}_q$ is computed by an operation $\widehat{\bfa}_q
= B^{-1} \bfa_q$ known as \FTRAN. Update of the dual variables
$\widehat{\bfc}\subN^T$ (\UPDATEDUAL) is achieved using $\theta_q$ and
$\widehat{\bfa}_p^T$. The update of the DSE weights is given by
\[
\begin{array}{rll}
w_p & := w_p / \widehat{a}_{pq}^2 \\\noalign{\smallskip}
w_i & := w_i - 2 (\widehat{a}_{iq} / \widehat{a}_{pq}) \tau_i + (\widehat{a}_{iq} / \widehat{a}_{pq})^2 w_p &i\ne p
\end{array}
\]
This requires both the \FTRAN result $\widehat{\bfa}_q$ and the
solution of $\bftau = B^{-1} \widehat{\bfe}_p$. The latter is obtained
by another \FTRAN type operation, known as \FTRANDSE.

Following a BFRT ratio test, if $\SetF$ is not empty, then all the
variables with indices in $\SetF$ are flipped, and the primal basic solution
$\bfx\subB$ is further updated (another \UPDATEPRIMAL) by the result
of the \FTRANBFRT operation $\widehat{\bfa}\subF = B^{-1} \bfa\subF$,
where $\bfa\subF$ is a linear combination of the constraint columns
for the variables in $\SetF$.

\subsubsection{Scope for parallelisation}
The computational components identified above are summarised in Table~\ref{table:components}. This also gives the
average contribution to solution time for the LP test set used in
Section~\ref{sect:Results}.

\begin{table}[hbt]
\tbl{Major components of the dual revised simplex method and their percentage of overall solution time\label{table:components}}{%
\begin{tabular}{l | l | r}
\hline
\TableOneForceLineHeight
Components    & Brief description & Percentage  \\
\hline
\TableOneForceLineHeight
\INVERT       & Recompute $B^{-1}$&13.3\\
\nullbox{0pt}{6pt}
\UPDATEFACTOR & Update basis inverse $B_k^{-1}$ to $B_{k+1}^{-1}$ &2.3\\
\hline
\CHUZR        & Choose leaving variable $p$ &2.9\\
\hline 
\TableOneForceLineHeight
\BTRAN        & Solve for $\widehat{\bfe}_p^T = \bfe_p^T B^{-1}$ &8.7\\
\TableOneForceLineHeight
\PRICE        & Compute $\widehat{\bfa}_p^T = \widehat{\bfa}_p^T A$&18.4\\
\CHUZCONE     & Collect valid ratio test candidates&7.3\\
\CHUZCTWO     & Search for entering variable $p$&1.5\\
\hline
\TableOneForceLineHeight
\FTRAN        & Solve for $\widehat{\bfa}_q = B^{-1} \bfa_q$&10.8\\
\FTRANBFRT    & Solve for $\widehat{\bfa}\subF = B^{-1} \bfa\subF$&3.5\\
\nullbox{0pt}{6pt}
\FTRANDSE     & Solve for $\bftau = B^{-1} \widehat{\bfe}_p$ &26.4\\
\hline
\TableOneForceLineHeight
\UPDATEDUAL   & Update $\widehat{\bfc}^T$ using $\widehat{\bfa}_p^T$&\\
\UPDATEPRIMAL & Update $\bfx\subB$ using $\widehat{\bfa}_q$ or $\widehat{\bfa}\subF$ &4.8\\
\UPDATEWEIGHT & Update DSE weight using $\widehat{\bfa}_q$ and $\bftau$ &\\
\hline
\end{tabular}}
\end{table}

There is immediate scope for data parallelisation within \CHUZR,
\PRICE, \CHUZC and most of the update operations since they require
independent operations for each (nonzero) component of a vector.
Exploiting such parallelisation in \PRICE and \CHUZC has been reported by Bixby and
Martin~\cite{Bixby2000} who achieve speedup on a small group of LP
problems with relatively expensive \PRICE operations.  The scope for
task parallelism by overlapping \FTRAN and \FTRANDSE was considered by
Bixby and Martin but rejected as being disadvantageous
computationally.

\subsection{Dual suboptimization}
\label{background:subopt}

Suboptimization is one of the oldest variants of the revised simplex
method and consists of a major-minor iteration scheme. Within the
primal revised simplex method, suboptimization performs minor
iterations of the standard primal simplex method using small subsets
of columns from the reduced coefficient matrix
$\widehat{A}=B^{-1}A$. Suboptimization for the dual simplex method was
first set out by Rosander~\cite{Rosander1975} but no practical
implementation has been reported. It performs minor operations of the
standard dual simplex method, applied to small subsets of rows from
$\widehat{A}$.

\begin{enumerate}
  \item {\it Major optimality test}.  Choose index set $\SetP \subseteq
    \SetB$ of primal infeasible basic variables as potential leaving
    variables.  If no such indices can be chosen, the LP problem has
    been solved to optimality.
  \item {\it Minor initialisation.}  For each $p \in \SetP$, compute
    $\widehat{\bfe}_p^T = \bfe_p^T B^{-1}$.
  \item {\it Minor iterations.}
  \begin{enumerate}
    \item {\it Minor optimality test.} Choose and remove 
          a primal infeasible variable $p$ from $\SetP$. 
          If no such variable can be chosen, 
          the minor iterations are terminated. 
    \item {\it Minor ratio test.} As in the regular ratio test,
      compute $\widehat{\bfa}_p^T = \widehat{\bfe}_p^T A$ (\PRICE) then
      identify an entering variable $q$.
    \item {\it Minor update.} Update primal variables for the
      remaining candidates in set $\SetP$ only ($\bfx_\SetP$) and
      update all dual variables $\widehat{\bfc}\subN$.
  \end{enumerate}
  \item {\it Major update}.
  For the pivotal sequence identified during the minor iterations, 
  update the primal basic variables, DSE weights and representation of $B^{-1}$.
\end{enumerate}

Originally, suboptimization was proposed as a pivoting scheme
with the aim of achieving better pivot choices and advantageous
data affinity.  In modern revised simplex implementations, the DSE and
BFRT are together regarded as the best pivotal rules and the idea of
suboptimization has been largely forgotten.

However, in terms of parallelisation, suboptimization is attractive
because it provides more scope for parallelisation. For the primal
simplex algorithm, suboptimization underpinned the work of Hall and
McKinnon~\cite{Hall1996,Hall1998}.  For dual suboptimization the major
initialisation requires $s$ \BTRAN operations, where $s=|\SetP|$. Following
$t\le s$ minor iterations, the major
update requires $t$ \FTRAN operations, $t$ \FTRANDSE operations and up
to $t$ \FTRANBFRT operations. The detailed design of the
parallelisation scheme based on suboptimization is discussed in
Section~\ref{sect:PAMI}.

\subsection{Simplex update techniques}
\label{background:update}

Updating the basis inverse $B_k^{-1} $ to $ B_{k+1}^{-1}$ after the
basis change $B_{k+1} = B_k + (\bfa_q - B\bfe_p)\bfe_p^T$ is a crucial
component of revised simplex method implementations.  The standard
choices are the
relatively simple product form (PF) update~\cite{Orchard-Hays1968} or
the efficient Forrest-Tomlin (FT) update~\cite{Forrest1972}.  A
comprehensive report on simplex update techniques is given by Elble
and Sahinidis~\cite{ElSa12b} and novel techniques, some motivated by
the design and development of \pami, are described by Huangfu and
Hall~\cite{HuHa12}.  For the purpose of this report, the features
of all relevant update methods are summarised as follows.

\begin{itemize}
\item The {\it product form\/} (PF) update uses the \FTRAN result
$\widehat{\bfa}_q$, yielding $B_{k+1}^{-1} = E^{-1}B_k^{-1}$, where the inverse of $E
= I + (\widehat{\bfa}_q - \bfe_p) \bfe^T_p$, is readily
available.

\item The {\it Forrest-Tomlin\/} (FT) update assumes $B_k=L_kU_k$ and uses
both the partial \FTRAN result $\tilde{\bfa}_q = L_k^{-1}\bfa_q$ and
partial \BTRAN result $\tilde{\bfe}_p^T = {\bfe}_p^T U_k^{-1}$ to
modify 
$U_k$ and augment
$L_k$.

\item The {\it alternate product form\/} (APF) update~\cite{HuHa12} uses
the \BTRAN result $\widehat{\bfe}_p^T$ so that $B_{k+1}^{-1} =
B_k^{-1}T^{-1}$, where $T = I+(\bfa_q - \bfa_{p'})\widehat{\bfe}_p^T$
and $\bfa_{p'}$ is column $p$ of $B$. Again, $T$ is readily inverted.

\item Following suboptimization, the {\it collective Forrest-Tomlin\/} (CFT)
update~\cite{HuHa12} updates $B_k^{-1}$ to $B_{k+t}^{-1}$
directly, using partial results obtained with $B_k^{-1}$ which are
required for simplex iterations. 
\end{itemize}

Although the direct update of the basis inverse from $B_k^{-1}$ to
$B_{k+t}^{-1}$ can be achieved easily via the PF or APF update, in
terms of efficiency for future simplex iterations, the collective FT
update is preferred to the PF and APF updates. The value of the APF update
within \pami is indicated in Section~\ref{sect:PAMI}.

\section{Parallelism across multiple iterations}\label{sect:PAMI}

This section introduces the design and implementation of the parallel
dual simplex scheme, \pami. It extends the suboptimization scheme of
Rosander~\cite{Rosander1975}, incorporating (serial) algorithmic
techniques and exploiting parallelism across multiple iterations.

The fundamental design of \pami was introduced by Hall and
Huangfu~\cite{Hall2012}, where it was referred to as ParISS. This
prototype implementation was based on the PF update and was relatively unsophisticated, both
algorithmically and computationally. Subsequent revisions
and refinements, incorporating the advanced algorithmic techniques
outlined in Section~\ref{sect:Background} as well as FT updates and
some novel features described in this section, have yielded a very much
more sophisticated and efficient implementation.

Section~\ref{pami:overview} provides an overview of the
parallelisation scheme of \pami and Section~\ref{subsect:pami-ftran}
details the task parallel \FTRAN operations in the major update stage
and how to simplify it. A novel
candidate quality control scheme for the minor optimality test is discussed in Section~\ref{pami:chuzr}.

\subsection{Overview of the \pami framework}
\label{pami:overview}

This section details the general \pami parallelisation scheme with
reference to the suboptimization framework introduced in
Section~\ref{background:subopt}.  

\subsubsection{Major optimality test}

The major optimality test involves only major \CHUZR operations in
which $s$ candidates are chosen (if possible) using the DSE framework. In
\pami the value of $s$ is the number of processors being used. It is a
vector-based operation which can be easily parallelised,
although its overall computational cost is not significant since it is
only performed once per major operation. However, the algorithmic design of
\CHUZR is important and Section~\ref{pami:chuzr} discusses it in
detail.

\subsubsection{Minor initialisation}

The minor initialisation step computes the \BTRAN results for (up to
$s$) potential candidates to leave the basis.  This is the first of
the task parallelisation opportunities provided by the suboptimization
framework.

\subsubsection{Minor iterations}
\label{pami:minor_iterations}

There are three main operations in the minor iterations.

\begin{itemize}
\item[(a)] Minor \CHUZR simply chooses the best candidates from the set
$\SetP$. Since this is computationally trivial, exploitation of
parallelism is not considered. However, consideration must be given to
the likelihood that the attractiveness of the best remaining candidate
in $\SetP$ has dropped significantly. In such circumstances, it may not be desirable to allow this
variable to leave the basis. This consideration leads to a candidate
quality control scheme introduced in Section~\ref{pami:chuzr}.

\item[(b)] The minor ratio test is a major source of parallelisation and
performance improvement.  Since the \BTRAN result is known (see
below), the minor ratio test consists of \PRICE, \CHUZCONE and
\CHUZCTWO.  The \PRICE operation is a sparse matrix-vector product and
\CHUZCONE is a one-pass selection based on the result of \PRICE.  In
the actual implementation, they can share one parallel initialisation.
On the other hand, \CHUZCTWO often involves multiple iterations of
recursive selection which, if exploiting parallelism, requires many
synchronisation operations.  According to the component profiling in
Table~\ref{table:components}, \CHUZCTWO is a relative cheap operation
thus, in \pami, it is not parallelised.
Data parallelism is exploited in \PRICE and \CHUZCONE by partitioning
the variables across the processors before any simplex iterations are
performed. This is done randomly with the aim of achieving load
balance in \PRICE.

\item[(c)] The minor update consists of the update of dual
variables and the update of \BTRAN results. The former is performed in
the minor update because the dual variables are required in the ratio
test of the next minor iteration. It is simply a vector addition and
represents immediate data parallelism.  The updated \BTRAN result
$\bfe_i^TB_{k+1}^{-1}$ is obtained by observing that it is given by
the APF update as
$\bfe_i^TB_k^{-1}T^{-1}=\widehat{\bfe}_i^TT^{-1}$. Exploiting the
structure of $T^{-1}$ yields a vector operation which may be
parallelised. After the \BTRAN results have been updated, the DSE
weights of the remaining candidates are recomputed directly at little cost.
\end{itemize}

\subsubsection{Major update}

Following $t$ minor iterations, the major update step concludes the
major iteration.  It consists of three types of operation: up to
$3t$ \FTRAN operations (including \FTRANDSE and \FTRANBFRT), the vector-based
update of primal variables and DSE weights, and
update of the basis inverse representation.

The number of \FTRAN operations cannot be fixed {\em a priori\/} since
it depends on the number of minor iterations and the number involving a
non-trivial BFRT. A simplification of the group of \FTRAN{}s is
introduced in~\ref{subsect:pami-ftran}.

The updates of all primal variables and DSE weights (given the particular
vector $\bftau = B^{-1} \widehat{\bfe}_p$) are vector-based data
parallel operations.

The update of the invertible representation of $B$ is performed using
the collective FT update unless it is desirable or necessary to
perform \INVERT to reinvert $B$. Note that both of these operations
are performed serially. Although the (collective) FT update is
relatively cheap (see Table~\ref{table:components}), so has little impact on performance, there is
significant processor idleness during the serial \INVERT.

\subsection{Parallelising three groups of \FTRAN operations}
\label{subsect:pami-ftran}

%
%
Within \pami, the pivot sequence $\{p_i, q_i\}_{i=0}^{t-1}$
identified in minor iterations yields up to $3t$ forward linear systems (where $t\leq
s$).
%
%
Computationally, there are three groups of \FTRAN operations, being $t$ regular
\FTRAN{}s for obtaining updated tableau columns $\widehat{\bfa}_q =
B^{-1} \bfa_q$ associated with the entering variable identified during
minor iterations; $t$ additional \FTRANDSE operations to obtain the
DSE update vector $\bftau = B^{-1} \widehat{\bfe}_p$ and 
\FTRANBFRT calculations to update the primal solution resulting from
bound flips identified in the BFRT.  Each system in a group is
associated with a different basis matrix, $B_k, B_{k+1}, \ldots,
B_{k+t-1}$.  For example the $t$ regular forward systems for obtaining
updated tableau columns are $\widehat{\bfa}_{q_0} =
B_k^{-1}{\bfa}_{q_0}, \widehat{\bfa}_{q_1} = B_{k+1}^{-1}{\bfa}_{q_1},
\ldots, \widehat{\bfa}_{q_{t-1}} = B_{k+t-1}^{-1}{\bfa}_{q_{t-1}}$.

%
%
For the regular \FTRAN and \FTRANDSE operations, the $i^\mathrm{th}$
linear system (which requires $B_{k+i}^{-1}$) in each group, is solved
by applying $B_k^{-1}$ followed by $i-1$ PF transformations given by
$\widehat{\bfa}_{q_j}, j < i$ to bring the result up to date.  The
operations with $B_k^{-1}$ and PF transformations are referred to as
the inverse and update parts respectively.  The multiple inverse parts
are easily arranged as a task parallel computation. The update part of
the regular \FTRAN operations requires results of other forward
systems in the same group and thus cannot be performed as task
parallel calculations. However, it is possible and valuable to exploit
data parallelism when applying individual PF updates when $\widehat{\bfa}_{q_i}$
is large and dense.  For the \FTRANDSE group it is possible to exploit
task parallelism fully if this group of computations is performed
after the regular \FTRAN.  However, when implementing \pami, both
\FTRANDSE and regular \FTRAN are performed together to increase the
number of independent inverse parts in the interests of load
balancing.

%
%
The group of up to $t$ linear systems associated with BFRT is slightly
different from the other two groups of systems.  Firstly,
there may be anything between none and $t$ linear systems depending
how many minor iterations are associated with actual bound flips.
More importantly, the results are only used to update the
values of the primal variables $\bfx\subB$ by simple vector addition.
This can be expressed as a single operation
\begin{equation}
\bfx\subB := \bfx\subB + 
\sum_{i=0}^{t-1} B_{k+i}^{-1} {\bfa}_{\scriptscriptstyle {Fi}}
=\bfx\subB + \sum_{i=0}^{t-1} 
\left(
\prod_{j=i-1}^{0} E_j^{-1} B_{k}^{-1} {\bfa}_{\scriptscriptstyle {Fi}}
\right)
\label{equ:pami:bfrt:pf}
\end{equation}
where one or more of ${\bfa}_{\scriptscriptstyle {Fi}}$ may be a zero
vector. If implemented using the regular PF update, each \FTRANBFRT operation starts
from the same basis inverse $B_k^{-1}$ but finishes with different
numbers of PF update operations.  Although these operations are
closely related, they cannot be combined. However, if the APF update is
used, so $B_{k+i}^{-1}$ can be expressed as
\[
B_{k+i}^{-1} =  B_k^{-1} T_0^{-1} \ldots T_{i-1}^{-1},
\]
the primal update equation~(\ref{equ:pami:bfrt:pf}) can be
rewritten as
\begin{equation}
\bfx\subB := 
\bfx\subB + \sum_{i=0}^{t-1} 
\left(
B_{k}^{-1} \prod_{j=0}^{i-1} T_j^{-1} {\bfa}_{\scriptscriptstyle {Fi}}
\right)
=
\bfx\subB + B_{k}^{-1} 
\left(
\sum_{i=0}^{t-1} 
\prod_{j=0}^{i-1} T_j^{-1} {\bfa}_{\scriptscriptstyle {Fi}}
\right)
\label{equ:pami:bfrt:apf}
\end{equation}
where the $t$ linear systems start with a cheap APF update part and
finish with a {\em single\/} $B_k^{-1}$ operation applied to the
combined result.  This approach greatly reduces the total serial cost
of solving the forward linear systems associated with BFRT. An
additional benefit of this combination is that the \UPDATEPRIMAL
operation is also reduced to a single operation after the combined
\FTRANBFRT.

%
%
By combining several potential \FTRANBFRT operations into one, the number of
forward linear systems to be solved is reduced to $2t + 1$, or $2t$
when no bound flips are performed.  An additional benefit of this
reduction is that, when $t\le s-1$, the total number of forward linear
systems to be solved is less than $2s$, so that each
of the $s$ processors will solve at most two linear systems.  However, when $t
= s$ and \FTRANBFRT is nontrivial, one of the $s$ processors is
required to solve three linear systems, while the other processors are
assigned only two, resulting in an ``orphan task''. To avoid this
situation, the number of minor iterations is limited to $t = s-1$ if
bound flips have been performed in the previous $s-2$ iterations.

\begin{figure}[hbt]
\centering
\begin{tikzpicture}
\draw node at (0,0)  [block, minimum height=1.8cm] {FTRAN \\ BFRT};
\draw node at (3cm,0)  [block, minimum height=1.8cm] {FTRAN \\[2ex] FTRAN \\ DSE};
\draw node at (6cm,0)  [block, minimum height=1.8cm] {FTRAN \\[2ex] FTRAN \\ DSE};
\draw node at (9cm,0)  [block, minimum height=1.8cm] {FTRAN \\[2ex] FTRAN \\ DSE};
\draw [dashed](2cm, 0.2cm) -- (4cm, 0.2cm);
\draw [dashed](5cm, 0.2cm) -- (7cm, 0.2cm);
\draw [dashed](8cm, 0.2cm) -- (10cm, 0.2cm);
\draw node at (4.5cm,-2cm)  [block, text width=11cm] {\phantom{U}FTRAN \hspace{3ex} UPDATE};
\draw [gray, dashed](1.5cm, 1cm) -- (1.5cm, -3cm);
\draw [gray, dashed](4.5cm, 1cm) -- (4.5cm, -3cm);
\draw [gray, dashed](7.5cm, 1cm) -- (7.5cm, -3cm);
\end{tikzpicture}

\caption{Task parallel scheme of all \FTRAN operations in \pami}
\label{figure:ftranparallel}
\end{figure}
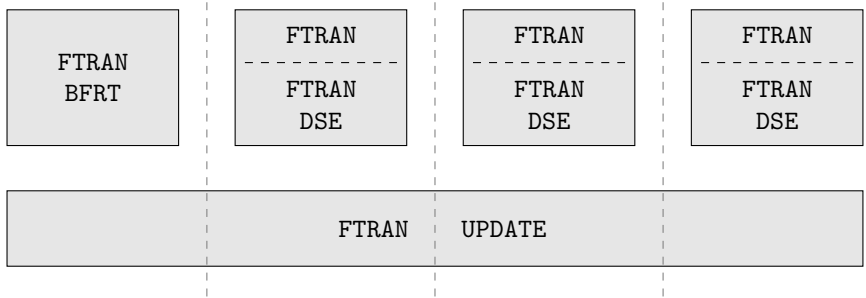

%
%
The arrangement of the task parallel \FTRAN operations discussed above
is illustrated in Figure~\ref{figure:ftranparallel}.  In the actual
implementation, the $2t + 1$ \FTRAN operations are all started
the same time as parallel tasks, and the processors are left to decide
which ones to perform.

\subsection{Candidate persistence and quality control in \CHUZR}
\label{pami:chuzr}

%
%
Major \CHUZR forms the set $\SetP$ and minor \CHUZR chooses
candidates from it. The design of \CHUZR contributes significantly to the
serial efficiency of suboptimization schemes so merits careful
discussion.

%
%
When suboptimization is performed, the candidate chosen to leave the
basis in the first minor iteration is the same as would have been
chosen without suboptimization. Thereafter, the candidates remaining
in $\SetP$ may be less attractive than the most attractive of the
candidates not in $\SetP$ due to the former becoming less attractive
and/or the latter becoming more attractive. Indeed, some candidates in
$\SetP$ may become unattractive. If candidates in the original $\SetP$
do not enter the basis then the work of their \BTRAN operations (and
any subsequent updates) is wasted. However, if minor iterations choose
less attractive candidates to leave the basis the number of simplex
iterations required to solve a given LP problem can be expected to
increase. Addressing this issue of {\em candidate persistence\/} is
the key algorithmic challenge when implementing suboptimization. The
number of candidates in the initial set $\SetP$ must be decided, and a
strategy determined for assessing whether a particular candidate should
remain in $\SetP$.

For load balancing during the minor initialisation, the initial number
of candidates $s=|\SetP|$ should be an integer multiple of the number
of processors used. Multiples larger than one yield better load
balance due to the greater amount of work to be parallelised,
particularly before and after the minor iterations, but practical experience with \pami prototypes demonstrated clearly that this is more
than offset by the amount of wasted computation and an increase in the number of
iterations required to solve the problem. Thus, for \pami, $s$ was
chosen to be eight, whatever the number of processors. 

%
%
During minor iterations, after updating the primal activities of the
variables given by the current set $\SetP$, the attractiveness of 
$\alpha_p$ for each $p\in\SetP$ is assessed relative to its initial value $\alpha_p^i$ by
means of a {\it cutoff\/} factor $\psi>0$. Specifically, if
\[
\alpha_p < \psi \alpha_p^i,
\]
then index $p$ is removed from $\SetP$. Clearly if the variable
becomes feasible or unattractive ($\alpha_p\le0$) then it is dropped whatever the
value of $\psi$.

To determine the value of $\psi$ to use in \pami, a series of
experiments was carried out using a reference set of 30 LP problems given in
Table~\ref{table:lpset} of Section~\ref{subsect:TestProblems}, with cutoff ratios ranging from 1.001 to 0.01.
Computational results are presented in Table~\ref{table:cutoffexp}
which gives the (geometric) mean speedup factor and the number of
problems for which the speedup factor is respectively 1.6, 1.8 and
2.0.

\begin{table}[hbt]
\tbl{Experiments with different cutoff factor for controlling candidate quality in \pami\label{table:cutoffexp}}{%
\begin{tabular}{l | rrrr}
\hline
 cutoff ($\psi$) &  speedup & $\# 1.6$ speedup &$\# 1.8$ speedup & $\# 2.0$ speedup\\
\hline
   1.001 &     1.12 &        1 &        1 &        0\\
   0.999 &     1.52 &       11 &        7 &        5\\
    0.99 &     1.54 &       13 &        6 &        4\\
    0.98 &     1.53 &       15 &        8 &        5\\
    0.97 &     1.48 &       11 &        6 &        5\\
    0.96 &     1.52 &       12 &        8 &        6\\
    0.95 &     1.49 &       13 &        8 &        4\\
    0.94 &     1.56 &       13 &        8 &        4\\
    0.93 &     1.47 &       13 &        9 &        4\\
    0.92 &     1.52 &       14 &        7 &        4\\
    0.91 &     1.52 &       14 &        5 &        3\\
     0.9 &     1.50 &       12 &        9 &        4\\
     0.8 &     1.46 &       13 &        9 &        3\\
     0.7 &     1.46 &       15 &        9 &        4\\
     0.6 &     1.44 &       11 &        8 &        6\\
     0.5 &     1.42 &       13 &        5 &        3\\
     0.2 &     1.36 &       10 &        6 &        4\\
     0.1 &     1.29 &       10 &        7 &        3\\
    0.05 &     1.16 &        9 &        4 &        2\\
    0.02 &     1.28 &       10 &        6 &        2\\
    0.01 &     1.22 &        8 &        5 &        3\\
\hline
\end{tabular}}
\end{table}

The cutoff ratio $\psi = 1.001$ corresponds to a special situation,
in which only candidates associated with improved attractiveness are
chosen. As might be expected, the speedup with this value of $\psi$ is
poor.
The cutoff ratio $\psi = 0.999$ corresponds to a boundary situation
where candidates whose attractiveness decreases are dropped. An
mean speedup of~1.52 is achieved.

For various cutoff ratios in the range $0.9 \leq \psi \leq 0.999$,
there is no really difference in the performance of \pami: the
mean speedup and larger speedup counts are relatively stable.
Starting from $\psi = 0.9$, decreasing the cutoff factor results in a
clear decrease in the mean speedup, although the larger speedup
counts remain stable until $\psi = 0.5$.

In summary, experiments suggest that any value in interval $[0.9,
  0.999]$ can be chosen as the cutoff ratio, with \pami using the
median value $\psi = 0.95$.

\subsection{Hyper-sparse LP problems}

In the discussions above, when exploiting data parallelism in vector
operations it is assumed that one independent scalar calculation must be performed for most of the components of the vector. For example, in \UPDATEDUAL and \UPDATEPRIMAL a multiple of the component is added to the corresponding component of another vector. In \CHUZR and \CHUZCONE the component (if nonzero) is used to compute and then compare a ratio. Since these scalar calculations need not be performed for zero components of the vector, when
the LP problem exhibits hyper-sparsity this is exploited by efficient serial
implementations~\cite{Hall2005}. When the cost of the serial vector operation
is reduced in this way it is no longer efficient to exploit data
parallelism so, when the density of the vector is below a certain
threshold, \pami reverts to serial computation. The performance of \pami is not sensitive to the thresholds of 5\%--10\% which are used.

\section{Single iteration parallelism}\label{sect:SIP}

This section introduces a relative simple approach to exploiting
parallelism within a single iteration of the dual revised simplex
method, yielding the parallel scheme \sip.  Our approach is a
significant development of the work of Bixby and
Martin~\cite{Bixby2000} who parallelised only the \PRICE, \CHUZC and
\UPDATEDUAL operations, having rejected the task parallelism of
\FTRAN and \FTRANDSE as being computationally disadvantageous.

Our serial simplex solver \hsol has an additional \FTRANBFRT component
for the bound-flipping ratio test. However, naively exploiting task
parallelism by simply overlapping this with \FTRAN and \FTRANDSE is
inefficient since the latter is seen in Table~\ref{table:components}
to be relatively expensive. This is due to the RHS of \FTRANDSE being
$\widehat{\bfe}_p$, which is dense relative to the RHS vectors
$\bfa_q$ of \FTRAN and $\bfa\subF$ of \FTRANBFRT. There is also no
guarantee in a particular iteration that \FTRANBFRT will be required.

\begin{figure}[hbt]
\centering
\begin{tikzpicture}[yscale=1]
\draw node at (0,0)  [block] (A) {CHUZR};
\draw node at (0,-2cm) [block] (B) {BTRAN};
\draw [->] (A) -- (B) node [pos=0.5, right] {$p$};
\draw node at (0,-6cm) [block, minimum height=5cm] (C) {FTRAN\\ DSE \\[2ex] $(\bftau = B^{-1} \widehat{\bfe}_p)$};
\draw node at (3cm,-4cm) [block] (D0){CHUZC1 \\ (Logical)};
\draw node at (7.5cm,-4cm) [block, text width=5cm] (D){SPMV + CHUZC1 \\ (Structural)};
\draw [->] (B) -- (C) node [pos=0.5, right] {$\widehat{\bfe}_p$};
\draw [->] (B) -- (D0) node [pos=0.5, right,  xshift=2ex] {$\widehat{\bfe}_p$};
\draw [->] (B) -- (D) node [pos=0.5, right, xshift=4ex] {$\widehat{\bfe}_p$};
\draw node at (3cm,-6cm) [block] (E) {CHUZC2};
\draw [->] (D0) -- (E) node [pos=0.5, right] {$\SetJ_{(L)}$};
\draw [->] (D) -- (E) node [pos=0.5, right,  xshift=2ex] {$\SetJ_{(S)}$};
\draw node at (3cm,-8cm) [block] (F) {FTRAN};
\draw [->] (E) -- (F) node [pos=0.5, right] {$q$}; 
\draw node at (6cm,-8cm) [block] (G) {FTRAN \\ BFRT};
\draw [->, dashed] (E) -- (G) node [pos=0.5, right,xshift=1ex] {$\cal F$};
\draw node at (9cm,-8cm) [block] (H) {UPDATE\\ DUAL};
\draw [->] (E) -- (H) node [pos=0.5, right, xshift=2ex] {$\theta_d$};
\draw node at (0,-10cm) [block]  (I){UPDATE \\ WEIGHT};
\draw [->] (C) -- (I) node [pos=0.5, right] {$\bftau$};
\draw [->] (F) -- (I) node [pos=0.8, above] {$\widehat{\bfa}_q$};
\draw node at (3,-10cm) [block] (J) {UPDATE \\ PRIMAL};
\draw [->] (F) -- (J) node [pos=0.5, right] {$\widehat{\bfa}_q$};
\draw [->, dashed] (G) -- (J) node [pos=0.8, above] {$\widehat{\bfa}\subF$};
\draw [gray, dashed](1.5cm, 1cm) -- (1.5cm, -11cm);
\draw [gray, dashed](4.5cm, 1cm) -- (4.5cm, -11cm);
\draw [gray, dashed](7.5cm, 1cm) -- (7.5cm, -11cm);
\end{tikzpicture}

\caption{\sip data dependency and parallelisation scheme}
\label{figure:sip}
\end{figure}
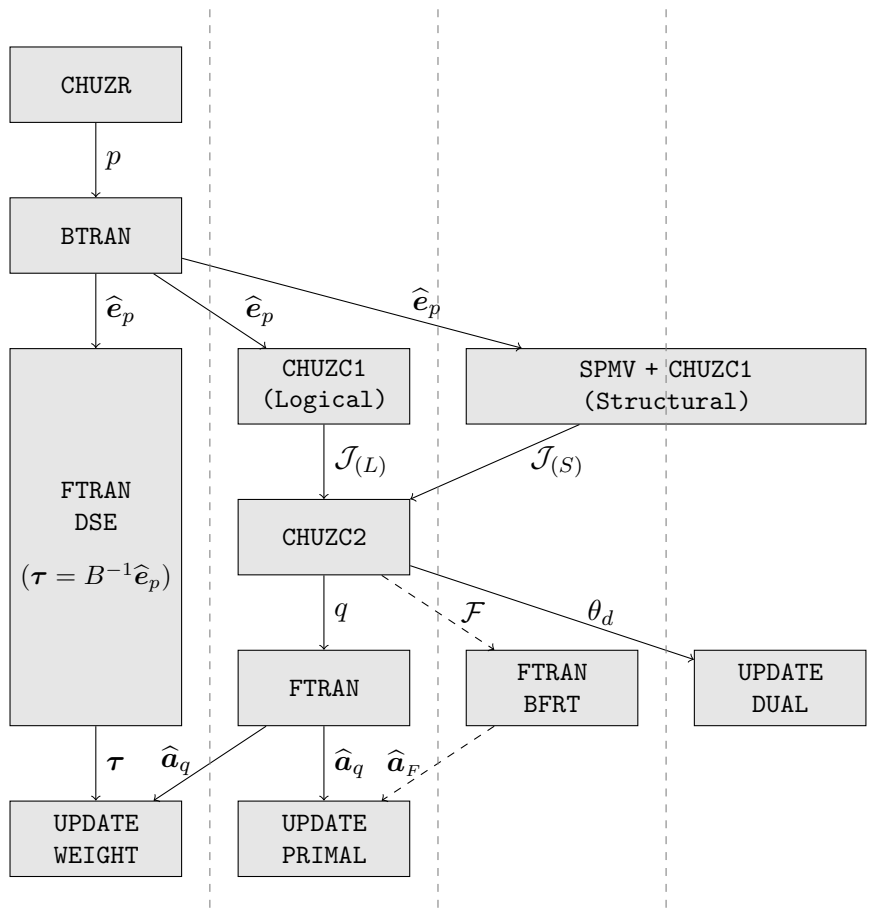

The mixed parallelisation scheme of \sip is illustrated in
Figure~\ref{figure:sip}, which also indicates the data dependency for
each computational component. Note that during \CHUZCONE there is a
distinction between the operations for the original (structural) variables and those for the logical (slack) variables,
since the latter correspond to an identity matrix in $A$. 
Thereafter, one processor performs \FTRAN in parallel with (any)
\FTRANBFRT on another processor and \UPDATEDUAL on a third. The scheme assumes
at least four processors but with more than four only the parallelism
in \PRICE and \CHUZC is enhanced.

\section{Computational results}\label{sect:Results}

\subsection{Test problems}
\label{subsect:TestProblems}
%
%
Throughout this report, the performance of the simplex solvers is
assessed using a reference set of 30 LP problems. Most of
these are taken from a comprehensive list of 
representative LP problems~\cite{MittelmannBenchmarksURL} maintained
by Mittelmann.

\begin{table}[h!]
\tbl{The reference set of 30 LP problems with hyper-sparsity measures\label{table:lpset}}{%
\begin{tabular}{l |rrr|rr}
\hline
\LP{Model           } &    \#row &    \#col &    \#nnz &  \FTRAN &   \BTRAN \\
\hline
\LP{cre-b           } &     9648 &    72447 &   256095 &     100 &       83 \\
\LP{dano3mip\_lp    } &     3202 &    13873 &    79655 &       1 &        6 \\
\LP{dbic1           } &    43200 &   183235 &  1038761 &     100 &       83 \\
\LP{dcp2            } &    32388 &    21087 &   559390 &     100 &       97 \\
\LP{dfl001          } &     6071 &    12230 &    35632 &      34 &       57 \\
\LP{fome12          } &    24284 &    48920 &   142528 &      45 &       58 \\
\LP{fome13          } &    48568 &    97840 &   285056 &     100 &       98 \\
\LP{ken-18          } &   105127 &   154699 &   358171 &     100 &      100 \\
\LP{l30             } &     2701 &    15380 &    51169 &      10 &        8 \\
\LP{Linf\_520c      } &    93326 &    69004 &   566193 &      10 &       11 \\
\LP{lp22            } &     2958 &    13434 &    65560 &      13 &       22 \\
\LP{maros-r7        } &     3136 &     9408 &   144848 &       5 &       13 \\
\LP{mod2            } &    35664 &    31728 &   198250 &      46 &       68 \\
\LP{ns1688926       } &    32768 &    16587 &  1712128 &      72 &      100 \\
\LP{nug12           } &     3192 &     8856 &    38304 &       1 &       20 \\
\LP{pds-40          } &    66844 &   212859 &   462128 &     100 &       98 \\
\LP{pds-80          } &   129181 &   426278 &   919524 &     100 &       99 \\
\LP{pds-100         } &   156243 &   505360 &  1086785 &     100 &       99 \\
\LP{pilot87         } &     2030 &     4883 &    73152 &      10 &       19 \\
\LP{qap12           } &     3192 &     8856 &    38304 &       2 &       15 \\
\LP{self            } &      960 &     7364 &  1148845 &       0 &        2 \\
\LP{sgpf5y6         } &   246077 &   308634 &   828070 &     100 &      100 \\
\LP{stat96v4        } &     3174 &    62212 &   490473 &      73 &       31 \\
\LP{stormG2-125     } &    66185 &   157496 &   418321 &     100 &      100 \\
\LP{stormG2-1000    } &   528185 &  1259121 &  3341696 &     100 &      100 \\
\LP{stp3d           } &   159488 &   204880 &   662128 &      95 &       70 \\
\LP{truss           } &     1000 &     8806 &    27836 &      37 &        2 \\
\LP{watson\_1       } &   201155 &   383927 &  1052028 &     100 &      100 \\
\LP{watson\_2       } &   352013 &   671861 &  1841028 &     100 &      100 \\
\LP{world           } &    35510 &    32734 &   198793 &      41 &       61 \\
\hline
\end{tabular}}
\end{table}

The problems in this reference set reflect the wide spread of LP
properties and revised simplex characteristics, including the
dimension of the linear systems (number of rows), the density of the
coefficient matrix (average number of non-zeros per column), and the
extent to which they exhibit hyper-sparsity (indicated by the last two
columns).  These columns, headed \FTRAN and \BTRAN, give the
proportion of the results of \FTRAN and \BTRAN with a density below
$10\%$, the criterion used to measure hyper-sparsity by Hall and
McKinnon~\cite{Hall2005} who consider an LP problem to be hyper-sparse
if the occurrence of such hyper-sparse results is greater than
$60\%$. According to this measurement, half of the reference set are
hyper-sparse. Since all problems are sparse, it is convenient to use
the term ``dense'' to refer to those which are not hyper-sparse.

The performance of \pami and \sip is assessed using experiments
performed on a workstation with 16 (Intel Xeon E5620, 2.4GHz)
processors, using eight for the parallel calculations. Numerical
results are given in Tables~\ref{table:TimeIters}
and~\ref{table:speedup}, where mean values of speedup or other relative
performance measures are computed geometrically. The relative
performance of solvers is also well illustrated using the performance
profiles in Figures~\ref{figure:HsolPamiClp}--\ref{figure:XpressCplex}.

\begin{sidewaystable}
\tbl{Iteration time (ms) and computational component profiling (the percentage of overall solution time) when solving LP problems with \hsol\label{table:compproffull}}{%
\begin{tabular}{lr|rrr|rr|rrrr|rr}
\hline
\LP{Model           } & Iter. Time & \CHUZR &\CHUZCONE &\CHUZCTWO &  \PRICE  &  \UPDATE &   \BTRAN &   \FTRAN &{\sc f-dse}&{\sc f-bfrt} & \INVERT & \OTHER \\
\hline
\LP{cre-b           } &      565 &      0.8 &     20.1 &      4.4 &     42.9 &      6.9 &      4.7 &      1.7 &     11.3 &      1.5 &      4.3 &      1.4 \\
\LP{dano3mip\_lp    } &      885 &      1.8 &     21.2 &      3.0 &     35.5 &      5.3 &      6.4 &      6.9 &     11.7 &      0.3 &      6.2 &      1.7 \\
\LP{dbic1           } &     2209 &      0.5 &     22.5 &      3.1 &     33.6 &      5.8 &      5.7 &      6.5 &     14.8 &      3.2 &      3.1 &      1.2 \\
\LP{dcp2            } &      509 &      6.5 &      3.9 &      1.7 &      8.7 &      7.3 &      5.4 &     18.1 &     28.4 &     10.4 &      7.4 &      2.2 \\
\LP{dfl001          } &      595 &      4.1 &      8.1 &      1.0 &     17.9 &     11.2 &     10.8 &     13.0 &     20.7 &      6.2 &      5.2 &      1.8 \\
\LP{fome12          } &      971 &      7.9 &      5.1 &      0.6 &     12.4 &      6.8 &     12.3 &     14.5 &     24.0 &      7.1 &      7.9 &      1.4 \\
\LP{fome13          } &     1225 &     10.1 &      4.2 &      0.5 &     10.6 &      5.6 &     11.4 &     13.5 &     26.4 &      6.7 &      9.6 &      1.4 \\
\LP{ken-18          } &      126 &      5.3 &      2.9 &      0.6 &      5.2 &      2.2 &      7.9 &     11.0 &     24.4 &      3.8 &     32.4 &      4.3 \\
\LP{l30             } &     1081 &      0.8 &     14.1 &      9.9 &     24.0 &      6.3 &      8.6 &      9.0 &     12.9 &      4.1 &      8.5 &      1.8 \\
\LP{Linf\_520c      } &    26168 &      1.5 &      2.3 &      0.1 &     11.8 &      4.0 &     16.6 &     19.7 &     23.2 &      0.0 &     19.2 &      1.6 \\
\LP{lp22            } &      888 &      2.0 &     10.9 &      2.0 &     23.3 &      8.4 &      9.4 &     10.4 &     14.9 &      6.8 &     10.0 &      1.9 \\
\LP{maros-r7        } &     1890 &      0.8 &      2.8 &      0.2 &     10.2 &      2.7 &     17.5 &     15.3 &     20.6 &      0.0 &     27.4 &      2.5 \\
\LP{mod2            } &     1214 &      4.2 &      7.5 &      1.0 &      9.9 &      8.5 &     11.5 &     17.4 &     29.1 &      5.4 &      4.0 &      1.5 \\
\LP{ns1688926       } &     1806 &      2.0 &      0.1 &      0.0 &      2.9 &      4.8 &      3.3 &     31.4 &     44.1 &      0.0 &      6.5 &      4.9 \\
\LP{nug12           } &     1157 &      1.6 &      7.4 &      1.1 &     16.3 &      6.9 &     11.6 &     12.4 &     16.7 &      5.8 &     18.1 &      2.1 \\
\LP{pds-40          } &      302 &      3.4 &      7.5 &      1.9 &     19.2 &      5.1 &     10.8 &     10.3 &     23.2 &      4.4 &     12.0 &      2.2 \\
\LP{pds-80          } &      337 &      3.7 &      6.6 &      1.8 &     19.8 &      3.9 &     10.5 &      9.1 &     23.7 &      3.9 &     15.0 &      2.0 \\
\LP{pds-100         } &      360 &      3.5 &      7.0 &      1.8 &     18.6 &      3.7 &     10.4 &      9.0 &     24.1 &      3.8 &     16.0 &      2.1 \\
\LP{pilot87         } &      918 &      1.2 &      5.1 &      0.8 &     17.9 &      4.4 &     12.0 &     12.9 &     17.4 &      7.6 &     17.9 &      2.8 \\
\LP{qap12           } &     1229 &      1.5 &      7.5 &      1.0 &     16.2 &      6.6 &     12.1 &     12.3 &     16.7 &      5.9 &     18.4 &      1.8 \\
\LP{self            } &     8350 &      0.0 &      1.4 &      0.2 &     39.6 &      0.2 &      7.0 &      6.5 &      7.0 &      0.0 &     33.9 &      4.2 \\
\LP{sgpf5y6         } &      491 &      1.3 &      0.3 &      0.1 &      0.2 &      0.1 &      5.0 &      2.3 &     80.7 &      0.0 &      8.4 &      1.6 \\
\LP{stat96v4        } &     2160 &      0.4 &     12.4 &      4.9 &     67.6 &      1.7 &      2.4 &      1.7 &      4.3 &      0.6 &      2.2 &      1.8 \\
\LP{stormG2-125     } &      115 &      5.2 &      0.8 &      0.2 &      1.7 &      0.9 &      4.4 &      8.3 &     48.7 &      0.1 &     26.7 &      3.0 \\
\LP{stormG2-1000    } &      650 &      1.5 &      0.1 &      0.0 &      0.3 &      1.3 &      3.5 &      6.1 &     70.6 &      0.0 &     14.6 &      2.0 \\
\LP{stp3d           } &     4325 &      1.6 &     10.7 &      0.9 &     19.2 &      7.6 &     13.5 &     12.0 &     27.0 &      3.9 &      2.4 &      1.2 \\
\LP{truss           } &      415 &      1.1 &     17.1 &      2.0 &     53.8 &      5.0 &      5.0 &      3.7 &      7.1 &      0.0 &      3.5 &      1.7 \\
\LP{watson\_1       } &      210 &      4.3 &      0.7 &      0.2 &      1.0 &      1.2 &      5.7 &      6.0 &     54.4 &      3.5 &     19.6 &      3.4 \\
\LP{watson\_2       } &      161 &      5.5 &      0.3 &      0.0 &      0.4 &      0.8 &      4.6 &      7.7 &     35.2 &      5.0 &     34.5 &      6.0 \\
\LP{world           } &     1383 &      3.8 &      8.7 &      1.3 &     10.9 &      8.6 &     11.6 &     16.5 &     28.0 &      5.5 &      3.7 &      1.4 \\
\hline
\LP{Average         } &      867 &      2.9 &      7.3 &      1.5 &     18.4 &      4.8 &      8.7 &     10.8 &     26.4 &      3.5 &     13.3 &      2.3 \\
\hline
\end{tabular}}
\end{sidewaystable}

\begin{sidewaystable}
\tbl{Solution time and iteration counts for \hsol, \pami, \sip, \clp and \cplex\label{table:TimeIters}}{%
\begin{tabular}{l |rrrrrr|rrrrr}
\hline
& \multicolumn{6}{c|}{Solution time} &  \multicolumn{5}{c}{Iteration counts} \\
\hline
 \LP{Model}           & \hsol   &\pamione&\pamieight & \sip & \clp & \cplex & \hsol   & \pami & \sip & \clp & \cplex \\
\hline
\LP{cre-b           } & 4.62 & 3.82 & 2.37 & 3.78 & 12.78 & 1.44 & 11599 & 10641 & 11632 & 26734 & 10912\\
\LP{dano3mip\_lp    } & 38.21 & 55.86 & 17.47 & 22.93 & 43.92 & 10.64 & 60161 & 47774 & 62581 & 64773 & 27438\\
\LP{dbic1           } & 52.43 & 111.22 & 39.24 & 44.43 & 542.62 & 27.64 & 35884 & 36373 & 37909 & 330315 & 46685\\
\LP{dcp2            } & 9.34 & 11.18 & 6.07 & 7.77 & 23.78 & 3.93 & 25360 & 24844 & 25360 & 43305 & 24036\\
\LP{dfl001          } & 11.74 & 17.80 & 6.31 & 8.47 & 13.13 & 7.89 & 26322 & 23668 & 26417 & 26866 & 21534\\
\LP{fome12          } & 71.74 & 116.92 & 42.26 & 56.50 & 54.22 & 50.58 & 103005 & 97646 & 101406 & 95142 & 85492\\
\LP{fome13          } & 186.35 & 271.72 & 113.39 & 148.27 & 122.58 & 156.90 & 209722 & 193928 & 204705 & 189503 & 177456\\
\LP{ken-18          } & 10.23 & 12.34 & 8.49 & 12.85 & 14.91 & 5.37 & 107471 & 106646 & 107467 & 106812 & 81952\\
\LP{l30             } & 7.93 & 17.48 & 6.24 & 6.04 & 7.14 & 5.60 & 10290 & 11433 & 10389 & 8934 & 10793\\
\LP{Linf\_520c      } & 2329.49 & 6402.00 & 2514.32 & 1699.63 & 6869.00 & 11922.00 & 132244 & 127468 & 132244 & 226319 & 153027\\
\LP{lp22            } & 15.74 & 26.54 & 9.64 & 10.97 & 14.90 & 8.54 & 25080 & 25778 & 24888 & 22401 & 18474\\
\LP{maros-r7        } & 7.91 & 27.49 & 16.08 & 6.47 & 8.60 & 2.73 & 6025 & 6258 & 6025 & 5643 & 6585\\
\LP{mod2            } & 38.90 & 73.57 & 29.78 & 32.39 & 25.77 & 19.83 & 43386 & 43100 & 42944 & 39552 & 48134\\
\LP{ns1688926       } & 17.75 & 28.13 & 10.16 & 12.96 & 2802.23 & 15.38 & 13849 & 15455 & 13849 & 193565 & 7228\\
\LP{nug12           } & 88.37 & 142.20 & 50.05 & 76.70 & 288.70 & 58.61 & 108152 & 102429 & 118370 & 211658 & 92368\\
\LP{pds-40          } & 20.39 & 31.28 & 15.04 & 18.08 & 155.53 & 16.26 & 94914 & 92992 & 92888 & 147122 & 58578\\
\LP{pds-80          } & 46.54 & 85.58 & 39.57 & 45.01 & 583.12 & 39.58 & 197461 & 200694 & 195658 & 409923 & 124097\\
\LP{pds-100         } & 59.21 & 94.67 & 46.32 & 55.06 & 719.33 & 51.88 & 234184 & 231758 & 231570 & 554434 & 143383\\
\LP{pilot87         } & 4.93 & 7.92 & 3.28 & 3.73 & 5.66 & 5.61 & 7240 & 7390 & 7130 & 8918 & 12069\\
\LP{qap12           } & 111.93 & 123.70 & 43.46 & 134.40 & 168.50 & 58.43 & 128131 & 86418 & 205278 & 134570 & 90736\\
\LP{self            } & 28.02 & 47.44 & 22.35 & 16.28 & 20.43 & 29.07 & 4738 & 5429 & 4738 & 4659 & 12073\\
\LP{sgpf5y6         } & 111.75 & 153.94 & 53.18 & 174.71 & 188.91 & 5.00 & 348115 & 346042 & 347978 & 347526 & 59716\\
\LP{stat96v4        } & 101.35 & 161.92 & 44.24 & 51.10 & 131.66 & 50.62 & 72531 & 65440 & 72531 & 119002 & 87056\\
\LP{stormG2-125     } & 7.02 & 8.95 & 5.58 & 10.00 & 18.01 & 3.98 & 81869 & 82965 & 81869 & 92149 & 86526\\
\LP{stormG2-1000    } & 290.35 & 397.72 & 185.34 & 352.44 & 1018.35 & 105.44 & 658534 & 658338 & 658534 & 738319 & 783176\\
\LP{stp3d           } & 355.98 & 443.99 & 152.47 & 305.96 & 254.71 & 163.98 & 130689 & 97680 & 130276 & 126346 & 98914\\
\LP{truss           } & 5.69 & 7.93 & 3.24 & 3.63 & 3.68 & 2.80 & 18929 & 15987 & 18929 & 17561 & 19693\\
\LP{watson\_1       } & 35.70 & 43.89 & 25.82 & 47.30 & 133.56 & 21.34 & 238973 & 239301 & 239819 & 466774 & 208888\\
\LP{watson\_2       } & 37.96 & 44.21 & 26.95 & 50.65 & 1118.00 & 35.88 & 334733 & 331607 & 334494 & 498797 & 305197\\
\LP{world           } & 47.97 & 86.49 & 34.29 & 38.69 & 33.83 & 26.19 & 47104 & 44722 & 46742 & 46283 & 54656\\
\hline
\end{tabular}}
\end{sidewaystable}

\begin{table}
\tbl{Speedup of \pami and \sip with hyper-sparsity measures\label{table:speedup}}{%
\begin{tabular}{l |rrrr|rr}
\hline
&  \multicolumn{4}{c}{Speedup} & \multicolumn{2}{|c}{Hyper-sparsity}\\
\hline
\LP{Model}            & \texttt{p1/hsol}& \texttt{p8/p1}&  \texttt{p8/hsol} & \texttt{sip/hsol} & \FTRAN & \BTRAN \\
\hline
\LP{cre-b           } & 1.21 & 1.61 & 1.95 & 1.22 &    100 &       83 \\
\LP{dano3mip\_lp    } & 0.68 & 3.20 & 2.19 & 1.67 &      1 &        6 \\
\LP{dbic1           } & 0.47 & 2.83 & 1.34 & 1.18 &    100 &       83 \\
\LP{dcp2            } & 0.84 & 1.84 & 1.54 & 1.20 &    100 &       97 \\
\LP{dfl001          } & 0.66 & 2.82 & 1.86 & 1.39 &     34 &       57 \\
\LP{fome12          } & 0.61 & 2.77 & 1.70 & 1.27 &     45 &       58 \\
\LP{fome13          } & 0.69 & 2.40 & 1.64 & 1.26 &    100 &       98 \\
\LP{ken-18          } & 0.83 & 1.45 & 1.20 & 0.80 &    100 &      100 \\
\LP{l30             } & 0.45 & 2.80 & 1.27 & 1.31 &     10 &        8 \\
\LP{Linf\_520c      } & 0.36 & 2.55 & 0.93 & 1.37 &     10 &       11 \\
\LP{lp22            } & 0.59 & 2.75 & 1.63 & 1.43 &     13 &       22 \\
\LP{maros-r7        } & 0.29 & 1.71 & 0.49 & 1.22 &      5 &       13 \\
\LP{mod2            } & 0.53 & 2.47 & 1.31 & 1.20 &     46 &       68 \\
\LP{ns1688926       } & 0.63 & 2.77 & 1.75 & 1.37 &     72 &      100 \\
\LP{nug12           } & 0.62 & 2.84 & 1.77 & 1.15 &      1 &       20 \\
\LP{pds-40          } & 0.65 & 2.08 & 1.36 & 1.13 &    100 &       98 \\
\LP{pds-80          } & 0.54 & 2.16 & 1.18 & 1.03 &    100 &       99 \\
\LP{pds-100         } & 0.63 & 2.04 & 1.28 & 1.08 &    100 &       99 \\
\LP{pilot87         } & 0.62 & 2.41 & 1.50 & 1.32 &     10 &       19 \\
\LP{qap12           } & 0.90 & 2.85 & 2.58 & 0.83 &      2 &       15 \\
\LP{self            } & 0.59 & 2.12 & 1.25 & 1.72 &      0 &        2 \\
\LP{sgpf5y6         } & 0.73 & 2.89 & 2.10 & 0.64 &    100 &      100 \\
\LP{stat96v4        } & 0.63 & 3.66 & 2.29 & 1.98 &     73 &       31 \\
\LP{stormG2-125     } & 0.78 & 1.60 & 1.26 & 0.70 &    100 &      100 \\
\LP{stormG2-1000    } & 0.73 & 2.15 & 1.57 & 0.82 &    100 &      100 \\
\LP{stp3d           } & 0.80 & 2.91 & 2.33 & 1.16 &     95 &       70 \\
\LP{truss           } & 0.72 & 2.45 & 1.76 & 1.57 &     37 &        2 \\
\LP{watson\_1       } & 0.81 & 1.70 & 1.38 & 0.75 &    100 &      100 \\
\LP{watson\_2       } & 0.86 & 1.64 & 1.41 & 0.75 &    100 &      100 \\
\LP{world           } & 0.55 & 2.52 & 1.40 & 1.24 &     41 &       61 \\
\hline
\LP{Mean            } & 0.64 & 2.34 & 1.51 & 1.15 \\
\hline
\end{tabular}}
\end{table}

\subsection{Performance of \pami}

\begin{figure}[hbt]
\centering
\begin{tikzpicture}[xscale=1.2,yscale=3.6,font=\small]
%
%
\draw [thick                           ]
( 1.00, 0.10)--( 1.18, 0.10)--
( 1.18, 0.13)--( 1.20, 0.13)--
( 1.20, 0.17)--( 1.31, 0.17)--
( 1.31, 0.20)--( 1.32, 0.20)--
( 1.32, 0.23)--( 1.64, 0.23)--
( 1.64, 0.27)--( 2.23, 0.27)--
( 2.23, 0.30)--( 2.51, 0.30)--
( 2.51, 0.33)--( 2.63, 0.33)--
( 2.63, 0.37)--( 2.70, 0.37)--
( 2.70, 0.40)--( 3.43, 0.40)--
( 3.43, 0.43)--( 4.41, 0.43)--
( 4.41, 0.47)--( 5.38, 0.47)--
( 5.38, 0.50)--( 5.45, 0.50)--
( 5.45, 0.53)--( 6.01, 0.53)--
( 6.01, 0.57)--( 6.74, 0.57)--
( 6.74, 0.60)--( 7.47, 0.60)--
( 7.47, 0.63)--( 7.56, 0.63)--
( 7.56, 0.67)--(10.00, 0.67);
\draw [densely dotted                  ]
( 1.00, 0.07)--( 1.40, 0.07)--
( 1.40, 0.10)--( 1.46, 0.10)--
( 1.46, 0.13)--( 1.58, 0.13)--
( 1.58, 0.17)--( 1.61, 0.17)--
( 1.61, 0.20)--( 1.63, 0.20)--
( 1.63, 0.23)--( 1.76, 0.23)--
( 1.76, 0.27)--( 1.80, 0.27)--
( 1.80, 0.30)--( 1.84, 0.30)--
( 1.84, 0.33)--( 1.86, 0.33)--
( 1.86, 0.37)--( 1.92, 0.37)--
( 1.92, 0.40)--( 1.94, 0.40)--
( 1.94, 0.43)--( 2.13, 0.43)--
( 2.13, 0.47)--( 2.15, 0.47)--
( 2.15, 0.50)--( 2.21, 0.50)--
( 2.21, 0.53)--( 2.27, 0.53)--
( 2.27, 0.57)--( 2.42, 0.57)--
( 2.42, 0.60)--( 2.45, 0.60)--
( 2.45, 0.63)--( 2.57, 0.63)--
( 2.57, 0.67)--( 2.68, 0.67)--
( 2.68, 0.70)--( 2.70, 0.70)--
( 2.70, 0.73)--( 2.72, 0.73)--
( 2.72, 0.77)--( 2.94, 0.77)--
( 2.94, 0.80)--( 3.14, 0.80)--
( 3.14, 0.83)--( 3.48, 0.83)--
( 3.48, 0.87)--( 3.67, 0.87)--
( 3.67, 0.90)--( 3.90, 0.90)--
( 3.90, 0.93)--( 4.00, 0.93)--
( 4.00, 0.97)--( 4.54, 0.97)--
( 4.54, 1.00)--(10.00, 1.00);
\draw [thin                            ]
( 1.00, 0.00)--( 2.02, 0.00)--
( 2.02, 0.03)--( 2.36, 0.03)--
( 2.36, 0.07)--( 2.38, 0.07)--
( 2.38, 0.10)--( 2.44, 0.10)--
( 2.44, 0.13)--( 2.57, 0.13)--
( 2.57, 0.17)--( 2.89, 0.17)--
( 2.89, 0.20)--( 3.35, 0.20)--
( 3.35, 0.23)--( 3.43, 0.23)--
( 3.43, 0.27)--( 3.58, 0.27)--
( 3.58, 0.30)--( 3.62, 0.30)--
( 3.62, 0.33)--( 3.97, 0.33)--
( 3.97, 0.37)--( 4.14, 0.37)--
( 4.14, 0.40)--( 4.18, 0.40)--
( 4.18, 0.43)--( 4.26, 0.43)--
( 4.26, 0.47)--( 4.50, 0.47)--
( 4.50, 0.50)--( 4.93, 0.50)--
( 4.93, 0.53)--( 4.94, 0.53)--
( 4.94, 0.57)--( 4.98, 0.57)--
( 4.98, 0.60)--( 4.98, 0.60)--
( 4.98, 0.63)--( 5.05, 0.63)--
( 5.05, 0.67)--( 5.10, 0.67)--
( 5.10, 0.70)--( 5.13, 0.70)--
( 5.13, 0.73)--( 5.14, 0.73)--
( 5.14, 0.77)--( 5.15, 0.77)--
( 5.15, 0.80)--( 5.17, 0.80)--
( 5.17, 0.83)--( 5.26, 0.83)--
( 5.26, 0.87)--( 5.30, 0.87)--
( 5.30, 0.90)--( 5.94, 0.90)--
( 5.94, 0.93)--( 6.57, 0.93)--
( 6.57, 0.97)--( 6.99, 0.97)--
( 6.99, 1.00)--(10.00, 1.00);
\draw [thick, dotted                   ]
( 1.00, 0.83)--( 1.03, 0.83)--
( 1.03, 0.87)--( 1.18, 0.87)--
( 1.18, 0.90)--( 1.21, 0.90)--
( 1.21, 0.93)--( 1.35, 0.93)--
( 1.35, 0.97)--( 3.32, 0.97)--
( 3.32, 1.00)--(10.00, 1.00);
\draw (1,1.1)--(1,0)--(10.0,0)--(10.0,1.1);
\draw ( 1.00,0)--( 1.00,-0.02) node [below] { 1};
\draw ( 3.25,0)--( 3.25,-0.02) node [below] { 2};
\draw ( 5.50,0)--( 5.50,-0.02) node [below] { 3};
\draw ( 7.75,0)--( 7.75,-0.02) node [below] { 4};
\draw (10.00,0)--(10.00,-0.02) node [below] { 5};
\foreach \y in {0,20,40,60,80,100} \draw (1,\y/100)--(0.95,\y/100) node [left] {$\y$};
\foreach \y in {0,20,40,60,80,100} \draw (10.00,\y/100)--(10.05,\y/100);
%
%
\draw [thick                           ]( 1.00,-0.30)--( 1.50,-0.30) node [right,text=black] {Clp             };
\draw [densely dotted                  ]( 3.00,-0.30)--( 3.50,-0.30) node [right,text=black] {hsol            };
\draw [thin                            ]( 5.00,-0.30)--( 5.50,-0.30) node [right,text=black] {pami            };
\draw [thick, dotted                   ]( 7.00,-0.30)--( 7.50,-0.30) node [right,text=black] {pami8           };
\end{tikzpicture}

\caption{Performance profile of \clp, \hsol, \pami and \pamieight without preprocessing or crash}
\label{figure:HsolPamiClp}
\end{figure}
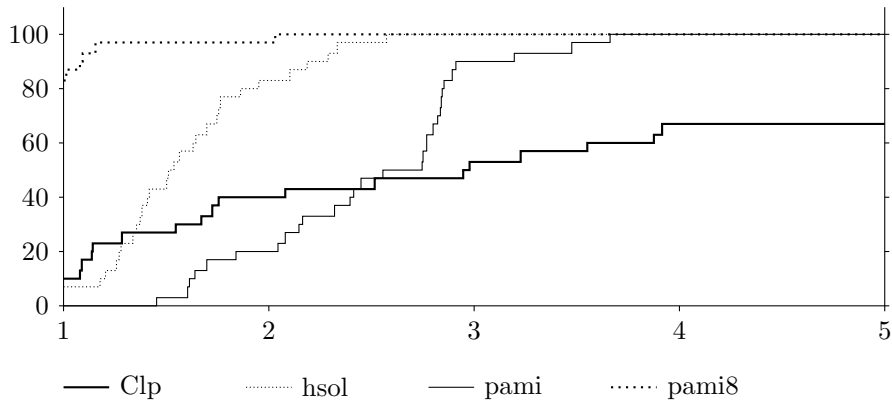

The efficiency of \pami is appropriately assessed in terms of
parallel speedup and performance relative to the sequential dual
simplex solver (\hsol) from which it was developed. The former
indicates the efficiency of the parallel implementation and the latter
measures the impact of suboptimization on serial performance. A high
degree of parallel efficiency would be of little value if it came at
the cost of severe serial inefficiency. The solution times for \hsol
and \pami running in serial, together with \pami running in parallel
with 8 cores, are listed in columns headed \hsol, \pamione and
\pamieight respectively in Table~\ref{table:TimeIters}. These results
are also illustrated via a performance profile in
Figure~\ref{figure:HsolPamiClp} which, to put the results in a broader
context, also includes \clp~1.15~\cite{ClpURL}, the world's leading
open-source solver. Note that since \hsol and \pami have no
preprocessing or crash facility, these are not used in the runs with
\clp.

The number of iterations required to solve a given LP problem can vary
significantly depending on the solver used and/or the algorithmic
variant used. Thus, using solution times as the sole measure of
computational efficiency is misleading if there is a significant
difference in iteration counts for algorithmic reasons. However, this
is not the case for \hsol and \pami. Observing that \pami identifies
the same sequence of basis changes whether it is run in serial or
parallel, relative to \hsol, the number of iterations required by
\pami is similar, with the mean relative iteration count of 0.96 being
marginally in favour of \pami. Individual relative iteration counts
lie in $[0.85, 1.15]$ with the exception of those for \LP{qap12},
\LP{stp3d} and \LP{dano3mip\_lp} which, being 0.67, 0.75 and 0.79
respectively, are significantly in favour of \pami. Thus, with the
candidate quality control scheme discussed in
Section~\ref{pami:chuzr}, suboptimization is seen not compromise the
number of iterations required to solve LP problems. Relative to \clp,
\hsol typically takes fewer iterations, with the mean relative
iteration count being 0.70 and extreme values of 0.07 for
\LP{ns1688926} and 0.11 for \LP{dbic1}. 

It is immediately clear from the performance profile in
Figure~\ref{figure:HsolPamiClp} that, when using 8 cores, \pami is superior
to \hsol which, in turn, is generally superior to \clp. Observe that
the superior performance of \pami on 8 cores relative to \hsol comes
despite \pami in serial being inferior to \hsol. Specifically, using
the mean relative solution times in Table~\ref{table:speedup}, \pami
on 8 cores is 1.51 times faster than \hsol, which is 2.29 times faster
than \clp. Even when taking into account that \hsol requires 0.70
times the iterations of \clp, the iteration speed of \hsol is seen to
be 1.60 times faster than \clp: \hsol is a high quality dual revised
simplex solver.

Since \hsol and \pami require very similar numbers of iterations, the
mean value of 0.64 for the inferiority of \pami relative to \hsol in
terms of solution time reflects the the lower iteration speed of \pami
due to wasted computation.
For more than $65\%$ of the reference set \pami is twice as fast in
parallel, with a mean speedup of 2.34. However, relative to \hsol,
some of this efficiency is lost due to overcoming the wasted
computation, lowering the mean relative solution time to 1.51.

For individual problems, there is considerable variance in the speedup
of \pami over \hsol, reflecting the variety of factors which affect
performance and the wide range of test problems.
For the two problems where \pami performs best in parallel, it is
flattered by requiring significantly fewer iterations than
\hsol. However, even if the speedups of 2.58 for \LP{qap12} and 2.33
for \LP{stp3d} are scaled by the relative iteration counts, the
resulting relative iteration speedups are still 1.74 and 1.75
respectively. However, for other problems where \pami performs well,
this is achieved with an iteration count which is similar to that of
\hsol. Thus the greater solution efficiency due to exploiting
parallelism is genuine.
Parallel \pami is not advantageous for all problems. Indeed, for
\LP{maros-r7} and \LP{Linf\_520c}, \pami is slower in parallel than
\hsol. For these two problems, serial \pami is slower than \hsol by
factors of 3.48 and 2.75 respectively. In addition, as can be seen in
Table~\ref{table:compproffull}, a significant proportion of the
computation time for \hsol is accounted for by \INVERT, which runs in
serial on one processor with no work overlapped.

Interestingly, there is no real relation between the performance of
\pami and problem hyper-sparsity: it shows almost same range of good,
fair and modest performance across both classes of problems,
although the more extreme performances are for dense problems. Amongst
hyper-sparse problems, the three where \pami performs best are
\LP{cre-b}, \LP{sgpf5y6} and \LP{stp3d}. This is due to the large
percentage of the solution time for \hsol accounted for by \PRICE
($42.9\%$ for \LP{cre-b} and $19.2\%$ for \LP{stp3d}) and \FTRANDSE
($80.7\%$ for \LP{sgpf5y6} and $27\%$ for \LP{stp3d}). In \pami, the
\PRICE and \FTRANDSE components can be performed efficiently as task
parallel and data parallel computations respectively, and therefore
the larger percentage of solution time accounted for by these
components yields a natural source of speedup.

\subsection{Performance of \sip}
\label{subsect:SipPerformance}

For \sip, the iteration counts are generally very similar to those of
\hsol, with the relative values lying in $[0.98, 1.06]$ except for the
two, highly degenerate problems \LP{nug12} and \LP{qap12} where \sip
requires 1.09 and 1.60 times as many iterations respectively. [Note
  that these two problems are essentially identical, differing only by
  row and column permutations.] It is clear from
Table~\ref{table:speedup} that the overall performance and mean
speedup (1.15) of \sip is inferior to that of \pami. This is because
\sip exploits only limited parallelism.

The worst cases when using \sip are associated with the hyper-sparse
LP problems where \sip typically results in a slowdown.
Such an example is \LP{sgpf5y6}, where the proportion of \FTRANDSE
is more than 80\% and the total proportion of \PRICE, \CHUZC, \FTRAN
and \UPDATEDUAL is less than 5\%.  Therefore, when performing
\FTRANDSE and the rest as task parallel operations, the overall
performance is not only limited by \FTRANDSE, but the competition for
memory access by the other components and the cost of setting up the
parallel environment will also slow down \FTRANDSE.

However, when applied to dense LP problems, the performance of \sip is
moderate and relatively stable.  This is especially so for those
instances where \pami exhibits a slowdown: for \LP{Linf\_520c},
\LP{maros-r7}, applying \sip achieves speedups of 1.31 and 1.12
respectively.

In summary, \sip, is a straightforward approach to parallelisation
which exploits purely single iteration parallelism and achieves
relatively poor speedup for general LP problems compared to \pami.
However, \sip is frequently complementary to \pami in achieving
speedup when \pami results in slowdown.

\subsection{Performance relative to \cplex and influence on \xpress}

\JHcomment{}{The overall performance of \pami is assessed via the performance
profile shown in Figure~\ref{figure:HsolPamiClp}.  As the figure
indicates, the performance of \pami is comparable with the dual
simplex implementation of \cplex~12.4~\cite{CplexURL}, a world-leading
commercial dual revised simplex solver, and clearly superior to that of
\clp~1.15~\cite{ClpURL}, the world's leading open-source solver.}

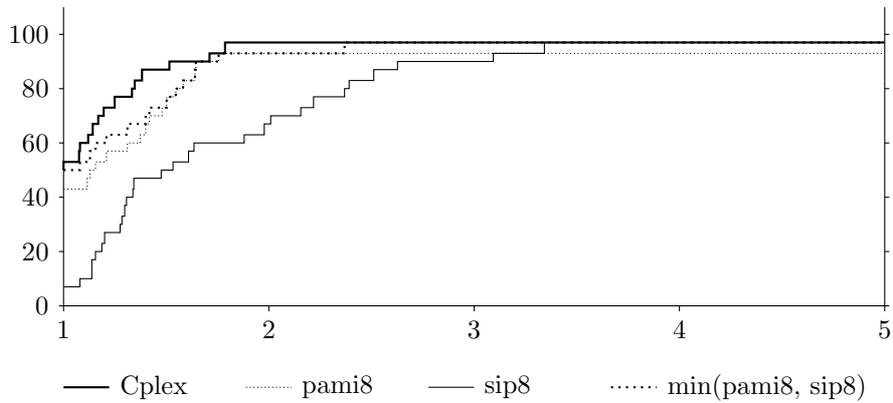
\begin{figure}[hbt]
\centering
\begin{tikzpicture}[xscale=1.2,yscale=3.6,font=\small]
%
%
\draw [thick                           ]
( 1.00, 0.50)--( 1.00, 0.50)--
( 1.00, 0.53)--( 1.17, 0.53)--
( 1.17, 0.57)--( 1.18, 0.57)--
( 1.18, 0.60)--( 1.27, 0.60)--
( 1.27, 0.63)--( 1.32, 0.63)--
( 1.32, 0.67)--( 1.38, 0.67)--
( 1.38, 0.70)--( 1.44, 0.70)--
( 1.44, 0.73)--( 1.56, 0.73)--
( 1.56, 0.77)--( 1.75, 0.77)--
( 1.75, 0.80)--( 1.78, 0.80)--
( 1.78, 0.83)--( 1.86, 0.83)--
( 1.86, 0.87)--( 2.16, 0.87)--
( 2.16, 0.90)--( 2.60, 0.90)--
( 2.60, 0.93)--( 2.77, 0.93)--
( 2.77, 0.97)--(10.00, 0.97);
\draw [densely dotted                  ]
( 1.00, 0.43)--( 1.26, 0.43)--
( 1.26, 0.47)--( 1.29, 0.47)--
( 1.29, 0.50)--( 1.35, 0.50)--
( 1.35, 0.53)--( 1.47, 0.53)--
( 1.47, 0.57)--( 1.70, 0.57)--
( 1.70, 0.60)--( 1.84, 0.60)--
( 1.84, 0.63)--( 1.90, 0.63)--
( 1.90, 0.67)--( 1.94, 0.67)--
( 1.94, 0.70)--( 2.08, 0.70)--
( 2.08, 0.73)--( 2.13, 0.73)--
( 2.13, 0.77)--( 2.23, 0.77)--
( 2.23, 0.80)--( 2.31, 0.80)--
( 2.31, 0.83)--( 2.44, 0.83)--
( 2.44, 0.87)--( 2.45, 0.87)--
( 2.45, 0.90)--( 2.70, 0.90)--
( 2.70, 0.93)--(10.00, 0.93);
\draw [thin                            ]
( 1.00, 0.07)--( 1.18, 0.07)--
( 1.18, 0.10)--( 1.31, 0.10)--
( 1.31, 0.13)--( 1.31, 0.13)--
( 1.31, 0.17)--( 1.35, 0.17)--
( 1.35, 0.20)--( 1.42, 0.20)--
( 1.42, 0.23)--( 1.45, 0.23)--
( 1.45, 0.27)--( 1.62, 0.27)--
( 1.62, 0.30)--( 1.64, 0.30)--
( 1.64, 0.33)--( 1.67, 0.33)--
( 1.67, 0.37)--( 1.69, 0.37)--
( 1.69, 0.40)--( 1.76, 0.40)--
( 1.76, 0.43)--( 1.77, 0.43)--
( 1.77, 0.47)--( 2.07, 0.47)--
( 2.07, 0.50)--( 2.20, 0.50)--
( 2.20, 0.53)--( 2.37, 0.53)--
( 2.37, 0.57)--( 2.43, 0.57)--
( 2.43, 0.60)--( 2.98, 0.60)--
( 2.98, 0.63)--( 3.20, 0.63)--
( 3.20, 0.67)--( 3.27, 0.67)--
( 3.27, 0.70)--( 3.60, 0.70)--
( 3.60, 0.73)--( 3.74, 0.73)--
( 3.74, 0.77)--( 4.08, 0.77)--
( 4.08, 0.80)--( 4.13, 0.80)--
( 4.13, 0.83)--( 4.40, 0.83)--
( 4.40, 0.87)--( 4.66, 0.87)--
( 4.66, 0.90)--( 5.71, 0.90)--
( 5.71, 0.93)--( 6.27, 0.93)--
( 6.27, 0.97)--(10.00, 0.97);
\draw [thick, dotted                   ]
( 1.00, 0.50)--( 1.18, 0.50)--
( 1.18, 0.53)--( 1.29, 0.53)--
( 1.29, 0.57)--( 1.35, 0.57)--
( 1.35, 0.60)--( 1.47, 0.60)--
( 1.47, 0.63)--( 1.70, 0.63)--
( 1.70, 0.67)--( 1.90, 0.67)--
( 1.90, 0.70)--( 1.94, 0.70)--
( 1.94, 0.73)--( 2.13, 0.73)--
( 2.13, 0.77)--( 2.23, 0.77)--
( 2.23, 0.80)--( 2.31, 0.80)--
( 2.31, 0.83)--( 2.44, 0.83)--
( 2.44, 0.87)--( 2.45, 0.87)--
( 2.45, 0.90)--( 2.70, 0.90)--
( 2.70, 0.93)--( 4.08, 0.93)--
( 4.08, 0.97)--(10.00, 0.97);
\draw (1,1.1)--(1,0)--(10.0,0)--(10.0,1.1);
\draw ( 1.00,0)--( 1.00,-0.02) node [below] { 1};
\draw ( 3.25,0)--( 3.25,-0.02) node [below] { 2};
\draw ( 5.50,0)--( 5.50,-0.02) node [below] { 3};
\draw ( 7.75,0)--( 7.75,-0.02) node [below] { 4};
\draw (10.00,0)--(10.00,-0.02) node [below] { 5};
\foreach \y in {0,20,40,60,80,100} \draw (1,\y/100)--(0.95,\y/100) node [left] {$\y$};
\foreach \y in {0,20,40,60,80,100} \draw (10.00,\y/100)--(10.05,\y/100);
%
%
\draw [thick                           ]( 1.00,-0.30)--( 1.50,-0.30) node [right,text=black] {Cplex           };
\draw [densely dotted                  ]( 3.00,-0.30)--( 3.50,-0.30) node [right,text=black] {pami8           };
\draw [thin                            ]( 5.00,-0.30)--( 5.50,-0.30) node [right,text=black] {sip8            };
\draw [thick, dotted                   ]( 7.00,-0.30)--( 7.50,-0.30) node [right,text=black] {min(pami8, sip8)};
\end{tikzpicture}

\caption{Performance profile of \cplex, \pamieight and \sipeight without preprocessing or crash}\label{figure:PamiSipCplex}
\end{figure}

Since commercial LP solvers are now highly developed it is, perhaps,
unreasonable to compare their performance with a research
code. However, this is done in Figure~\ref{figure:PamiSipCplex}, which
illustrates the performance of \cplex~12.4~\cite{CplexURL} relative to
\pamieight and \sipeight. Again, \cplex is run without preprocessing
or crash. Figure~\ref{figure:PamiSipCplex} also traces the performance
of the better of \pamieight and \sipeight, clearly illustrating that
\sip and \pami are frequently complementary in terms of achieving
speedup. Indeed, the performance of the better of \sip and \pami is
comparable with that of \cplex for the majority of the test
problems. For a research code this is a significant achievement.

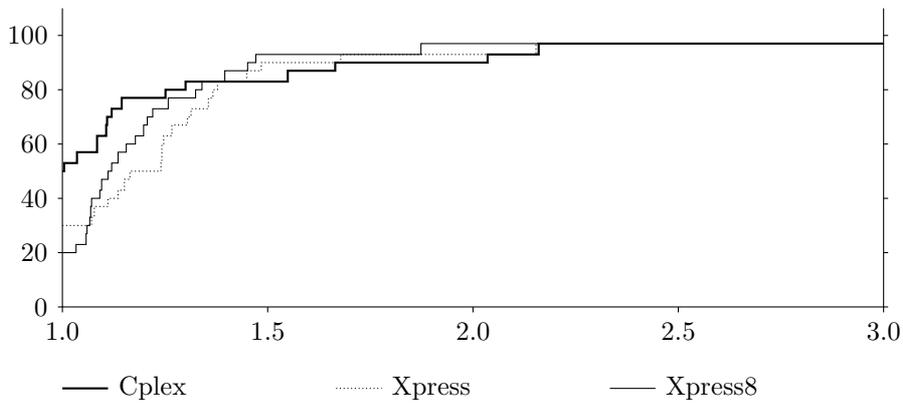
\begin{figure}[hbt]
\centering
\begin{tikzpicture}[xscale=1.2,yscale=3.6,font=\small]
%
%
\draw [thick                           ]
( 1.00, 0.50)--( 1.02, 0.50)--
( 1.02, 0.53)--( 1.16, 0.53)--
( 1.16, 0.57)--( 1.38, 0.57)--
( 1.38, 0.60)--( 1.38, 0.60)--
( 1.38, 0.63)--( 1.48, 0.63)--
( 1.48, 0.67)--( 1.49, 0.67)--
( 1.49, 0.70)--( 1.54, 0.70)--
( 1.54, 0.73)--( 1.65, 0.73)--
( 1.65, 0.77)--( 2.13, 0.77)--
( 2.13, 0.80)--( 2.35, 0.80)--
( 2.35, 0.83)--( 3.47, 0.83)--
( 3.47, 0.87)--( 3.99, 0.87)--
( 3.99, 0.90)--( 5.66, 0.90)--
( 5.66, 0.93)--( 6.22, 0.93)--
( 6.22, 0.97)--(10.00, 0.97);
\draw [densely dotted                  ]
( 1.00, 0.30)--( 1.32, 0.30)--
( 1.32, 0.33)--( 1.35, 0.33)--
( 1.35, 0.37)--( 1.50, 0.37)--
( 1.50, 0.40)--( 1.61, 0.40)--
( 1.61, 0.43)--( 1.68, 0.43)--
( 1.68, 0.47)--( 1.74, 0.47)--
( 1.74, 0.50)--( 2.08, 0.50)--
( 2.08, 0.53)--( 2.09, 0.53)--
( 2.09, 0.57)--( 2.09, 0.57)--
( 2.09, 0.60)--( 2.11, 0.60)--
( 2.11, 0.63)--( 2.20, 0.63)--
( 2.20, 0.67)--( 2.37, 0.67)--
( 2.37, 0.70)--( 2.41, 0.70)--
( 2.41, 0.73)--( 2.60, 0.73)--
( 2.60, 0.77)--( 2.65, 0.77)--
( 2.65, 0.80)--( 2.70, 0.80)--
( 2.70, 0.83)--( 3.02, 0.83)--
( 3.02, 0.87)--( 3.18, 0.87)--
( 3.18, 0.90)--( 4.05, 0.90)--
( 4.05, 0.93)--( 6.19, 0.93)--
( 6.19, 0.97)--(10.00, 0.97);
\draw [thin                            ]
( 1.00, 0.20)--( 1.15, 0.20)--
( 1.15, 0.23)--( 1.26, 0.23)--
( 1.26, 0.27)--( 1.27, 0.27)--
( 1.27, 0.30)--( 1.30, 0.30)--
( 1.30, 0.33)--( 1.31, 0.33)--
( 1.31, 0.37)--( 1.32, 0.37)--
( 1.32, 0.40)--( 1.41, 0.40)--
( 1.41, 0.43)--( 1.43, 0.43)--
( 1.43, 0.47)--( 1.50, 0.47)--
( 1.50, 0.50)--( 1.54, 0.50)--
( 1.54, 0.53)--( 1.61, 0.53)--
( 1.61, 0.57)--( 1.70, 0.57)--
( 1.70, 0.60)--( 1.80, 0.60)--
( 1.80, 0.63)--( 1.89, 0.63)--
( 1.89, 0.67)--( 1.93, 0.67)--
( 1.93, 0.70)--( 1.99, 0.70)--
( 1.99, 0.73)--( 2.16, 0.73)--
( 2.16, 0.77)--( 2.46, 0.77)--
( 2.46, 0.80)--( 2.53, 0.80)--
( 2.53, 0.83)--( 2.78, 0.83)--
( 2.78, 0.87)--( 3.03, 0.87)--
( 3.03, 0.90)--( 3.12, 0.90)--
( 3.12, 0.93)--( 4.93, 0.93)--
( 4.93, 0.97)--(10.00, 0.97);
\draw (1,1.1)--(1,0)--(10.0,0)--(10.0,1.1);
\draw ( 1.00,0)--( 1.00,-0.02) node [below] { 1.0};
\draw ( 3.25,0)--( 3.25,-0.02) node [below] { 1.5};
\draw ( 5.50,0)--( 5.50,-0.02) node [below] { 2.0};
\draw ( 7.75,0)--( 7.75,-0.02) node [below] { 2.5};
\draw (10.00,0)--(10.00,-0.02) node [below] { 3.0};
\foreach \y in {0,20,40,60,80,100} \draw (1,\y/100)--(0.95,\y/100) node [left] {$\y$};
\foreach \y in {0,20,40,60,80,100} \draw (10.00,\y/100)--(10.05,\y/100);
%
%
\draw [thick                           ]( 1.00,-0.30)--( 1.50,-0.30) node [right,text=black] {Cplex           };
\draw [densely dotted                  ]( 4.00,-0.30)--( 4.50,-0.30) node [right,text=black] {Xpress          };
\draw [thin                            ]( 7.00,-0.30)--( 7.50,-0.30) node [right,text=black] {Xpress8         };
\end{tikzpicture}

\caption{Performance profile of \cplex, \xpress and \xpresseight with preprocessing and crash}\label{figure:XpressCplex}
\end{figure}

Since developing and implementing the techniques described in this
paper, Huangfu has implemented them within the FICO \xpress simplex
solver, leading to FICO announcing via advertising copy and
blogs~\cite{QiBlogURL} that it has solved the long-standing problem of
parallelising the simplex method. The performance profile in
Figure~\ref{figure:XpressCplex} demonstrates that when it is
advantageous to run \xpress in parallel it enables FICO's solver to
match the serial performance of \cplex (which has no parallel simplex
facility). Note that for the results in in
Figure~\ref{figure:XpressCplex}, \xpress and \cplex were run with both
preprocessing and crash. The newly-competitive performance of parallel
\xpress relative to \cplex is also reflected in Mittelmann's
independent benchmarking~\cite{MittelmannBenchmarksURL}.

\section{Conclusions}\label{sect:Conclusions}

This report has introduced the design and development of two novel
parallel implementations of the dual revised simplex method. 

One relatively complicated parallel scheme (\pami) is based on a
less-known pivoting rule called suboptimization.  Although it provided
the scope for parallelism across multiple iterations, as a pivoting
rule suboptimization is generally inferior to the regular dual
steepest-edge algorithm.  Thus, to control the quality of the pivots,
which often declines during \pami, a {\it cutoff\/} factor is
necessary.  A suitable cutoff factor of 0.95, has been found via
series of experiments.  For the reference set \pami provides a mean
speedup of 1.51 which enables it to out-perform Clp, the best
open-source simplex solver.

The other scheme (\sip) exploits purely single iteration parallelism.
Although its mean speedup of 1.15 is worse than that of \pami, it is
frequently complementary to \pami in achieving speedup when \pami
results in slowdown.

Although the results in this paper are far from the linear speedup
which is the hallmark of many quality parallel implementations of
algorithms, to expect such results for an efficient implementation of
the revised simplex method applied to general large sparse LP problems
is unreasonable. The commercial value of efficient simplex
implementations is such that if such linear speedup were possible then
it would have been achieved years ago. A measure of the quality of the
\pami and \sip schemes discussed in this paper is that they have
formed the basis of refinements made by Huangfu to the Xpress solver
which have been considered noteworthy enough to be reported by FICO
and used in advertising copy. With the techniques described in this
paper, Huangfu has raised the performance of the Xpress parallel
revised simplex solver to that of the worlds best commercial simplex
solvers. In developing the first parallel revised simplex solver of
general utility and commercial importance, this work represents a
significant achievement in computational optimization.

\bibliographystyle{abbrv}
\bibliography{huangfu}

\end{document}